\documentclass[12pt]{amsart}
\usepackage{latexsym}
\usepackage{amssymb} 
\usepackage{mathrsfs}
\usepackage{amsmath}
\usepackage{latexsym}
\usepackage{delarray}
\usepackage{amssymb,amsmath,amsfonts,amsthm,
mathrsfs}

\usepackage{hyperref}
\renewcommand\eqref[1]{(\ref{#1})} 

\newtheorem{theorem}{Theorem}[section]

\newtheorem{proposition}[theorem]{Proposition}

\theoremstyle{definition}
\newtheorem{remark}[theorem]{Remark}

\newcommand{\wt}[1]{\widetilde{#1}}

\newcommand{\mb}[1]{\ensuremath{\mathbb{#1}}}
\newcommand{\N}{\mb{N}}

\newcommand{\R}{\mb{R}}

\renewcommand\N{{\mathbb N}_0}

\newcommand{\beq}{\begin{equation}}
\newcommand{\eeq}{\end{equation}}

\newcommand{\eps}{\varepsilon}

\newcommand{\lara}[1]{\langle #1 \rangle}

\renewcommand\N{{\mathbb N}_0}

\title[On higher order hyperbolic equations]
{On higher order hyperbolic equations with space-dependent coefficients: $C^\infty$ well-posedness and Levi conditions}

\author[Claudia Garetto]{Claudia Garetto}
\address{
  Claudia Garetto:
  \endgraf
School of Mathematical Sciences
  \endgraf
  Queen Mary University of London
  \endgraf
 Mile End Road, E1 4UJ, London
  \endgraf
  United Kingdom
  \endgraf
  {\it E-mail address} {\rm c.garetto@qmul.ac.uk}
  }

\thanks{The author was supported by the
EPSRC grant EP/V005529/2.
}
\date{}


\begin{document}

\maketitle

\begin{abstract}

This paper contributes to the wider study of hyperbolic equations with multiplicities. We focus here on some classes of higher order hyperbolic equations with space dependent coefficients in any space dimension. We prove Sobolev well-posedness of the corresponding Cauchy problem (with loss of derivatives due to the multiplicities) under suitable Levi conditions on the lower order terms. These conditions generalise the well known Olienik's conditions in \cite{O70} to orders higher than $2$. 
\end{abstract}

\section{Introduction}
The $C^\infty$ well-posedness of the Cauchy problem for hyperbolic equations with multiplicities has been a topic of great interested since the pioneering work of Oleinik in \cite{O70}. Note that the presence of multiplicities is often an obstacle to get $C^\infty$ well-posedness and differently from the strictly hyperbolic case, lower order terms play a relevant role in the analysis of these problems, see \cite{B, ColKi:02, CS, dAKi05, OT:84, PP} and references therein. The well-posedness result obtained by Oleinik holds for second order hyperbolic equations in variational form with smooth $(t,x)$-dependent coefficients and provides Sobolev well-posedness of any order with loss of derivatives. In detail, the Cauchy problem for a second order hyperbolic operator 
\begin{multline*}
Lu=u_{tt}-\sum_{i,j=1}^n (a_{ij}(t,x)u_{x_j})_{x_i}+\sum_{i=1}^n[(b_i(t,x)u_{x_i})_t+(b_i(t,x)u_t)_{x_i}]\\
+c(t,x)u_t+\sum_{i=1}^n d_i(t,x)u_{x_i}+e(t,x)u
\end{multline*}
with coefficients in $B^\infty([0,T]\times\R^n)$, the space of smooth functions with bounded derivatives of any order $k\in\N$, is $C^\infty$ well-posed if the lower order terms fulfil the following Oleinik's condition: there exist $A,C>0$ such that 
\[
\biggl[\sum_{i=1}^n d_i(t,x)\xi_i\biggr]^2\le C\biggl\{A\sum_{i,j=1}^na_{ij}(t,x)\xi_i\xi_j-\sum_{i,j=1}^n\partial_ta_{ij}(t,x)\xi_i\xi_j\biggr\},
\]
for all $t\in[0,T]$ and $x,\xi\in\R^n$. In the specific case of the wave operator
\[
\partial_t^2u-\sum_{i=1}^na_i(x)\partial^2_{x_i} u
\]
with $x$-dependent coefficients we have that 
\[
Lu=\partial_t^2u-\partial_{x_i}\biggl(\sum_{i=1}^na_i(x)\partial_{x_i} u\biggr)+\sum_{i=1}^n\partial_{x_i}a_i(x)\partial_{x_i}u
\]
and therefore Oleinik's condition is formulated as 
\beq
\label{est_OL}
\biggl[\sum_{i=1}^n \partial_{x_i}a_i(x)\xi_i\biggr]^2\le CA\sum_{i=1}^n a_i(x)\xi_i^2.
\eeq
Note that \eqref{est_OL} holds automatically by Glaeser's inequality if the coefficients $a_i$ are positive, at least of class $C^2$ and with bounded second order derivatives. Indeed,

{\bf Glaeser's inequality}: \emph{if $a\in C^2(\R^n)$, $a(x)\ge 0$ for all $x\in\R^n$ and 
\[
\sum_{i=1}^n\Vert \partial^2_{x_i}a\Vert_{L^\infty}\le M,
\]
for some constant $M>0$. Then,  
\[
|\partial_{x_i}a(x)|^2\le 2M a(x),
\]
for all $i=1,\dots,n$ and $x\in\R^n$.}
 
We can therefore state the following theorem.
\begin{theorem}
\label{theo_Oleinik}
The Cauchy problem 
\beq
\label{CP_wave}
\begin{split}
\partial_t^2u-\sum_{i=1}^na_i(x)\partial^2_{x_i} u&=f(t,x),\quad t\in[0,T], x\in\R^n,\\
u(0,x)&=g_0,\\
\partial_tu(0,x)&=g_1,
\end{split}
\eeq
where $a_i\in B^\infty(\R^n)$, $a_i\ge 0$, for all $i=1,\cdots, n$ and $f\in C([0,T], C^\infty_c(\R^n))$ is $C^\infty$ well-posed, i.e. given initial data $g_0,g_1\in C^\infty_c(\R^n)$ it has a unique smooth global solution on $[0,T]\times\R^n$.
\end{theorem}
It is not straightforward to extend Oleinik's result to higher order hyperbolic equations. This is due to technical difficulties arising from the higher number of roots and their multiplicities, so, to the best of our knowledge, the equivalent of Oleinik's condition for orders higher than 2 has not been formulated so far. However, mathematicians have investigated $C^\infty$ well-posedness for some special classes of equations: hyperbolic equations with $t$-dependent coefficients \cite{GarRuz:3, GarRuz:7, JT} and hyperbolic equations with coefficients in space dimension 1 \cite{ST, ST21}. Few results are also available for hyperbolic systems with multiplicities in diagonal \cite{KR2} and upper-triangular form \cite{GJR1, GJR2}. The general understanding of $C^\infty$ well-posedness for hyperbolic equations of any order with coefficients in $x\in\R^n$ is still open. In this paper we start to investigate $C^\infty$ well-posedness for higher order hyperbolic equations with coefficients in $x\in\R^n$. As remarked in \cite{ST21} well-posedness results holding in one space dimension do not necessarily hold in higher space dimension, however we prove here that when higher order hyperbolic equations are of a special form, namely without mixed $x$-derivatives in the principal part, then Levi conditions can be found for the lower order terms which guarantee $C^\infty$ well-posedness. These conditions can be regarded as an extension of Oleinik's conditions to orders higher than 2 and hold in any space dimension as well. 
For the third order equation 
\[
\partial_t^3u-\sum_{i=1}^n a_i(x)\partial_t\partial_{x_i}^2u\\
+\sum_{i=1}^n b_i(x)\partial_{x_i}^2 u+\sum_{i=1}^n  b_{2,i}(x)\partial_t\partial_{x_i}u+b_{3,n}(x)\partial_t^2u=f(t,x),
\]
our Levi conditions relate the coefficients $b_i$ and $b_{2,i}$ with $a_i$ and $\sqrt{a_i}$, respectively. Namely, $b_i=\lambda a_i$, for some $\lambda\in B^\infty(\R^n)$ and $|b_{2,i}|$ is bounded by $\sqrt{a_i}$, for all $i=1,\cdots, n$.
To explain our method, which employ ideas developed for the wave equation in \cite{G20}, we focus on special classes of hyperbolic equations of order $m=3$ and order $2m$ with $m\ge 2$ leaving the general treatment to a forthcoming paper which will employ pseudo-differential techniques rather than differential techniques.  
 
 The paper is organised as follows. In Section 2 we explain our method, based on reduction to a system of differential equations and construction of a symmetriser, on the wave equation toy model adding lower order terms to \eqref{CP_wave}. We show that the Levi conditions for lower order terms formulated by Oleinik can be also obtained from the system imposing that the matrix of the lower order terms is suitably estimated by the energy defined via the symmetriser. The extension of our method to third order equations is organised in two sections: Sections 3 and 4. We begin by analysing third order hyperbolic equations in space dimension 1 in Section 3 and we show that our method allows more general Levi conditions than the ones recently formulated in \cite{ST21}. We then pass to space dimension $n$ in Section 4. In Section 5,  we investigate a class of fourth order hyperbolic equations with $x$-dependent coefficients in $\R^n$ and a related class of equations of order $2m$ with $m\ge 2$.
Note that throughout the paper we work with real-valued functions and we look for real valued solutions. In all the Cauchy problems considered in this paper existence of the solution follows immediately from Nuij's approximation argument \cite{N68, ST21} and uniqueness is a direct consequence of the energy estimates. If the equation is of order $m$ we take lower order terms of order $m-1$ to perform a straightforward transformation into a system of differential equations however other lower order terms can be added by increasing the size of the system as in \cite{ST21}. A brief survey on the standard symmetriser for matrices in Sylvester form can be found in the appendix at the end of the paper.

\section{The case $m=2$}
For the sake of the reader we recall the method employed in \cite{G20} to prove the $C^\infty$ well-posedness of the Cauchy problem for the wave equation 
\[
\begin{split}
\partial_t^2u-\sum_{i=1}^na_i(x)\partial^2_{x_i} u&=f(t,x),\quad t\in[0,T], x\in\R^n,\\
u(0,x)&=g_0\in C^\infty_c(\R^n),\\
\partial_tu(0,x)&=g_1\in C^\infty_c(\R^n),
\end{split}
\]
where all the functions involved are real valued and $a_i(x)\ge 0$ for all $x\in\R^n$ and $i=1,\dots,n$.  Note that compactly supported initial data will enforce the solution $u$ to be compactly supported with respect to $x$ as well (finite propagation speed). We add to \eqref{CP_wave} lower order terms of any order. This leads to the equation
\[
\partial_t^2u-\sum_{i=1}^na_i(x)\partial^2_{x_i} u+\sum_{i=1}^nb_i(x)\partial_{x_i}u+c(x)\partial_t u+d(x)u=f(t,x),\quad t\in[0,T], x\in\R^n,
\]
that we can transform into a $(n+2)\times(n+2)$ system by setting
\[
U=\left(
	\begin{array}{c}
	U^{(0)}\\
	U^{(1)}
 \end{array}
	\right),
\]
where $U^{(0)}=u$ and
\[
U^{(1)}=(\partial_{x_1}u,\cdots,\partial_{x_n}u,\partial_t u)^T.
\]
In detail, we get
\[
\partial_tU=\sum_{i=1}^n A_i(x)\partial_{x_i}U+B(x)U+F,
\]
where the only non-zero entries of $A_i$ are $a_{1+i, n+2}=1$ and $a_{n+2, 1+i}=a_i$, for $i=1,\cdots,n$,
\[
B=\left(
	\begin{array}{ccccc}
	0 & 0 & \cdots & 0 &1\\
	0 & 0 & \cdots & 0 & 0\\
	\vdots & \vdots & \vdots & \vdots & \vdots\\
	0 & 0 & \cdots & 0 & 0\\
	-d &-b_1 & \cdots & -b_n & -c
  \end{array}
	\right)
\]
and 
\[
F=\left(
	\begin{array}{c}
	0\\
	0\\
	\vdots\\
	0\\
	f
  \end{array}
	\right).
\]
The initial data are given by $U(0)=(g_0,\partial_{x_1}g_0,\cdots, \partial_{x_n}g_0, g_1)^T$.

As an explanatory example we focus on the case $n=2$. We have
\[
A_1=\left(
	\begin{array}{cccc}
	0 & 0 & 0 & 0\\
	0& 0 & 0 & 1\\
	0& 0 & 0 & 0\\ 
	0 & a_1 & 0 & 0
	\end{array}
	\right),\quad 
A_2=\left(
	\begin{array}{cccc}
	0& 0 & 0 & 0\\
	0& 0 & 0 & 0\\ 
	0 & 0 & 0 & 1\\
	0 & 0 &  a_2 & 0 
	\end{array}
	\right).
\]
and
\[
B=\left(
	\begin{array}{cccc}
	0& 0 & 0 & 1\\
	0& 0 & 0 & 0\\ 
	0 & 0 & 0 & 0\\
	-d & -b_1 & -b_2  & -c 
	\end{array}
	\right).
\] 
Note that 
\[
Q=\left(
	\begin{array}{cccc}
	1& 0 & 0 & 0\\
	0& a_1 & 0 & 0\\ 
	0 & 0 & a_2 & 0\\
	0 & 0 & 0  & 1 
	\end{array}
	\right)
\]
is a symmetriser for both the matrices $A_1$ and $A_2$. Indeed,
\[
QA_1=\left(
	\begin{array}{cccc}
	0& 0 & 0 & 0\\
	0& 0 & 0 & a_1\\ 
	0 & 0 & 0 & 0\\
	0 & a_1 & 0  & 0
	\end{array}
	\right), \quad 
QA_2=\left(
	\begin{array}{cccc}
	0& 0 & 0 & 0\\
	0& 0 & 0 & 0\\ 
	0 & 0 & 0 & a_2\\
	0 & 0 & a_2  & 0
	\end{array}
	\right).
\]
Our system can be studied by using the energy $E=(QU,U)_{L^2}$ and employing the Glaeser's inequality as in \cite{G20}. In detail,
 \beq
 \label{first_deriv_E}
 \begin{split}
 &\frac{dE(t)}{dt}=(\partial_t(QU),U)_{L^2}+(QU,\partial_tU)_{L^2}\\
 &=-\sum_{k=1}^n(\partial_{x_k}(QA_k)U,U)_{L^2}+((QB+B^\ast Q)U,U)_{L^2}+2(QU,F)_{L^2}.
\end{split}
 \eeq
and therefore by analysing the term $((QB+B^\ast Q)U,U)_{L^2}$ we deduce how to formulate the Levi conditions on the lower order terms. We have
\[
QB+B^\ast Q=\left(
	\begin{array}{cccc}
	0& 0 & 0 & 1-d\\
	0& 0 & 0 & -b_1\\ 
	0 & 0 & 0 & -b_2\\
	1-d & -b_1 & -b_2  & -2c
	\end{array}
	\right)
\]
and by comparing\footnote{In the sequel given two functions $f=f(y)$ and $g=g(y)$ we use the notation $f\prec g$ if there exists a constant $C>0$  such that $f(y)\le C g(y)$, for all $y$.} $((QB+B^\ast Q)U,U)_{L^2}$ with $E(t)$ we get that
\[
\begin{split}
((QB+B^\ast Q)U,U)_{L^2}&= 2((1-d)U_1, U_4)-2(b_1U_2,U_4)-2(b_2U_3,U_4)-(2c U_4, U_4)\\
&{\prec} E(t)=(U_1,U_1)+(a_1 U_2, U_2)+(a_2 U_3, U_3)+(U_4, U_4)
\end{split}
\]
if 
\[
|d| \prec 1,\quad b_1^2 \prec a_1,\quad b_2^2\prec a_2,\quad |c| \prec 1.
 \]
These Levi conditions, which leads to Sobolev well-posedness with loss of derivatives, coincide with the well-known Oleinik's condition: 
\[
\biggl(\sum_{i=1}^n b_i(x)\xi_i\biggr)^2\le C \sum_{i=1}^n a_i(x)\xi_i^2, \qquad \text{for all $x,\xi\in\R^n$.}
\]
Indeed, in arbitrary space dimension $n$ we have 
\[
Q=\left(
	\begin{array}{cccccc}
	1& 0 & 0 & \cdots & 0 & 0\\
	0& a_1 & 0 & 0 & \cdots & 0\\ 
	0 & 0 & a_2 & 0 & \cdots & 0\\
	\vdots & \vdots & \vdots & \vdots & \vdots & \vdots\\
	0 & 0 & 0  &  \cdots & a_n & 0\\
	0 & 0 & 0  &  \cdots & \cdots & 1
	\end{array}
	\right)
\]
and 
\[
QB+B^\ast Q=\left(
	\begin{array}{cccccc}
	0& 0 & 0 & \cdots & 0 & 1-d\\
	0& 0 & 0 & \cdots & 0 &  -b_1\\ 
	0 & 0 & 0 & \cdots & 0 & -b_2\\
	\vdots & \vdots & \vdots & \vdots & \vdots & \vdots\\
	1-d & -b_1 & -b_2 & \cdots & -b_n & -2c
	\end{array}
	\right).
\]
which yields to
\[
|d| \prec 1,\quad b_1^2 \prec a_1,\quad b_2^2\prec a_2,\, \cdots, \, b_n^2\prec a_n,\quad |c| \prec 1.
\]
In the rest of the section, under the assumptions
\begin{itemize}
\item[(H)] \underline{the coefficients $a_i$ are non-negative and bounded for all $i=1,\dots, n$}\\
\underline{with bounded second order derivatives},
\item[(LC)] \underline{$b_i^2\prec a_i$ for all $i=1,\cdots,n$ and the lower order terms $c$ and $d$ are bounded},  
\end{itemize}
we analyse the terms in \eqref{first_deriv_E}. In few steps we will prove that our Cauchy problem is $C^\infty$ well-posed. This is the same conclusion reached by Oleinik however via a different analytical method which is more easily adaptable to higher order equations.\\[0.2cm]
{\bf Step 1: Estimate of the principal part}: by definition of the matrices $Q$ and $A_i$ we have that the only non-zero entries are the ones of indexes $i+1, n+2$ and $n+2, i+1$, respectively. They are both equal to $a_i$. So,
\[
((QA_i) U, U)_{L^2}=2(a_i U_{i+1}, U_{n+2})
\]
and 
\[
-\sum_{k=1}^n(\partial_{x_k}(QA_k)U,U)_{L^2}=-2\sum_{k=1}^n ((\partial_{x_k}a_k )U_{k+1}, U_{n+2}).
\]
It follows that 
\[
\biggl|-\sum_{k=1}^n(\partial_{x_k}(QA_k)U,U)_{L^2}\biggr|\le 2\sum_{k=1}^n |((\partial_{x_k}a_k )U_{k+1}, U_{n+2})|\le \sum_{k=1}^n\Vert \partial_{x_k}a_k U_{k+1}\Vert^2_{L^2} + n\Vert U_{n+2}\Vert^2_{L^2}.
\]
We now write $\Vert \partial_{x_k}a_kU_{k+1}\Vert^2_{L^2}$ as
\[
(\partial_{x_k}a_kU_{k+1},\partial_{x_k}a_{k}U_{k+1})_{L^2}=((\partial_{x_k}a_k)^2U_{k+1},U_{k+1})_{L^2}.
\]
Since $a_k\ge 0$ and $\sum_{j=1}^n\Vert \partial^2_{x_j}a_i\Vert_{L^\infty}\le M$ for all $i=1,\cdots,n$ by Glaeser's inequality  ($|\partial_{x_k}a_k(x)|^2\le 2Ma_k(x)$) we obtain the estimate
\[
 \Vert \partial_{x_k}a_kU_{k+1}\Vert^2_{L^2}\le 2M(a_kU_{k+1},U_{k+1})_{L^2}.
\]
Thus,
\[
\biggl|-\sum_{k=1}^n(\partial_{x_k}(QA_k)U,U)_{L^2}\biggr|\le 2M\sum_{k=1}^n(a_kU_{k+1},U_{k+1})_{L^2}+ n\Vert U_{n+2}\Vert^2_{L^2}\le \max(2M,n)E(t).
\]
\\[0.2cm]
{\bf Step 2: Estimate of the lower order terms}: from direct computations and by employing the Levi conditions (LC) we have
\[
\begin{split}
((QB+B^\ast Q)U,U)_{L^2}&=2((1-d)U_{n+2}, U_1)_{L^2}-2\sum_{k=1}^n (b_kU_{n+2},U_{k+1})_{L^2}-2(c U_{n+2}, U_{n+2})_{L^2}\\
&\prec \Vert U_1\Vert^2+\Vert U_{n+2}\Vert^2+\sum_{k=1}^n (a_k U_{k+1}, U_{k+1})+n\Vert U_{n+2}\Vert^2+\Vert U_{n+2}\Vert^2\\
&\prec E(t).
\end{split}
\] {\bf Step 3: Conclusion for $U_1$ and $U_{n+2}$}: there exists a constant $c'>0$ depending on $M$, the Levi conditions and the dimension $n$ such that 
\[
\frac{dE(t)}{dt}\le c' E(t)+\Vert f \Vert^2_{L^2}.
\]
By Gr\"onwall's lemma and the bound from below for the energy we obtain the following estimate for the entry $U_1=U^{(0)}=u$:  
\[
\begin{split}
\Vert u(t)\Vert_{L^2}^2=\Vert U_1\Vert^2_{L^2}&\le E(t)\le \biggl(E(0)+\int_{0}^t\Vert f(s)\Vert_{L^2}^2\, ds\biggr){\rm e}^{c't}\\
&\le  c''\biggl(\Vert g_0\Vert_{H^1}^2+\Vert g_1\Vert_{L^2}^2+\int_{0}^t\Vert f(s)\Vert_{L^2}^2\, ds\biggr).
\end{split}
\]
Analogously,  we have 
\[
\Vert \partial_tu(t)\Vert_{L^2}^2=\Vert U_{n+2}\Vert^2_{L^2}\le c''\biggl(\Vert g_0\Vert_{H^1}^2+\Vert g_1\Vert_{L^2}^2+\int_{0}^t\Vert f(s)\Vert_{L^2}^2\, ds\biggr).
\]
Note that in estimating $E(0)$ with the initial data we have used the fact that the coefficients $a_i$'s are bounded.\\[0.2cm]
{\bf Step 4: Estimates for $U_i$ with $i=2,\cdots, n+1$}. To get the well-posedness of the Cauchy problem
\[
\begin{split}
\partial_tU&=\sum_{i=1}^n A_i(x)\partial_{x_i}U+B(x)U+F,\\
U(0)&=(g_0,\partial_{x_1}g_0,\cdots, \partial_{x_n}g_0, g_1)^T
\end{split}
\]
we need to get an estimate for the remaining components of $U$, from $U_2$ to $U_{n+1}$. As in \cite{G20} we introduce $V=(\partial_{x_1}U, \cdots, \partial_{x_n}U)\in \R^{(n+2)n}$. Deriving with respect to $x$ we get
\[
\partial_t V=\sum_{i=1}^n \wt{A_i}(x)\partial_{x_i}V+ \wt{B}V+\wt{F},
\]
where $\wt{A_i}$ is a diagonal $(n+2)n\times (n+2)n$ matrix with $n$ repeated blocks $A_i$ on the diagonal, 
\[
{\wt{B}}=\left(
	\begin{array}{ccccc}
	\partial_{x_1}{A_1}+{B} & \partial_{x_1}{A_2} & \cdots & \cdots & \partial_{x_1}{A_n}\\
	\partial_{x_2}{A_1}& \partial_{x_2}{A_2}+{B} & \cdots & \cdots & \partial_{x_2}{A_n}\\
	\vdots & \vdots & \vdots & \vdots & \vdots\\
	\partial_{x_k}{A_1}&  \cdots &\partial_{x_k}{A_k}+{B} & \cdots & \partial_{x_k}{A_n}\\
	\vdots & \vdots & \vdots & \vdots & \vdots\\
	\partial_{x_n}{A_1}& \partial_{x_n}{A_2} & \cdots & \cdots & \partial_{x_n}{A_n}+{B}\\
	 
	\end{array}
	\right)
\]
and 
\[
{\wt{F}}=\nabla_x {F}+\left(
	\begin{array}{c}
	(\partial_{x_1}{B})U\\ 
	(\partial_{x_2}{B})U\\ 
	\vdots \\
	(\partial_{x_k}{B})U\\
	\vdots\\
	(\partial_{x_n}{B})U\\ 
	\end{array}
	\right)= \left(
	\begin{array}{c}
	\partial_{x_1}F\\ 
	\partial_{x_2}F\\ 
	\vdots \\
	\partial_{x_k}F\\
	\vdots\\
	\partial_{x_n}F\\ 
	\end{array}
	\right)+\left(
	\begin{array}{c}
	(\partial_{x_1}{B})U\\ 
	(\partial_{x_2}{B})U\\ 
	\vdots \\
	(\partial_{x_k}{B})U\\
	\vdots\\
	(\partial_{x_n}{B})U\\ 
	\end{array}
	\right).
\]
Arguing as in \cite{G20} Subsection 4.3, we make use of the energy $E(t)=(\wt{Q}V,V)_{L^2}$, where ${\wt{Q}}$ is a block-diagonal matrix with $n$ identical blocks equal to $Q$. In analogy with the system in $U$ we get 
\[
 \frac{dE(t)}{dt} =-\sum_{k=1}^n(\partial_{x_k}({\wt{Q}}{\wt{A_k}})V,V)_{L^2}+(({\wt{Q}}{\wt{B}}+{\wt{B}}^\ast{\wt{Q}})V,V)_{L^2}+2({\wt{Q}}V,{\wt{F}})_{L^2}.
\]
We now proceed with estimating this energy. Because of the block-diagonal structure of the matrices $\wt{A_k}$ and the symmetriser $\wt{Q}$ we argue as for the principal part of the system in $U$ and we get that 
\beq
\label{est_1_V}
\biggl|\sum_{k=1}^n(\partial_{x_k}({\wt{Q}}{\wt{A_k}})V,V)_{L^2}\bigg|\le c_1 E(t),
\eeq
for some constant $c_1>0$ depending on $M$ and the size of the matrices involved. This is clearly obtained under the hypothesis (H) and by applying the Glaeser's inequality. In order to estimate $(({\wt{Q}}{\wt{B}}+{\wt{B}}^\ast{\wt{Q}})V,V)_{L^2}$ it is sufficient to investigate the structure of $({\wt{Q}}{\wt{B}}V,V)_{L^2}$. This can be written as $(S_1V,V)_{L^2}+(S_2V,V)_{L^2}$, where  
\[
S_1=\wt{Q}\left(
	\begin{array}{ccccc}
	\partial_{x_1}{A_1} & \partial_{x_1}{A_2} & \cdots & \cdots & \partial_{x_1}{A_n}\\
	\partial_{x_2}{A_1}& \partial_{x_2}{A_2}& \cdots & \cdots & \partial_{x_2}{A_n}\\
	\vdots & \vdots & \vdots & \vdots & \vdots\\
	\partial_{x_k}{A_1}&  \cdots &\partial_{x_k}{A_k}& \cdots & \partial_{x_k}{A_n}\\
	\vdots & \vdots & \vdots & \vdots & \vdots\\
	\partial_{x_n}{A_1}& \partial_{x_n}{A_2} & \cdots & \cdots & \partial_{x_n}{A_n}\\
	 \end{array}
	\right)
\]
and $S_2$ is a block diagonal matrix with repeated blocks $QB$. The Levi conditions (LC) lead immediately to $(S_2V,V)_{L^2}\le c' E(t)$. It remains to estimate $(S_1 V,V)_{L^2}$. This means to deals with matrices of the type $Q\partial_{x_i}A_k$ for $i,k=1,\cdots,n$. By direct computations we easily see that $(S_1V,V)$ can be estimated blockwise with terms of the type, $(\partial_{x_i}a_kV_j, V_h)$ for a specific choice of indeces $j,h$ (see \cite{G20}, Subsection 4.3 for more details). For instance, when $n=2$ we get 
\[
\begin{split}
(S_1V,V)_{L^2}=& (\partial_{x_1}a_1 V_2,V_1)_{L^2}+(\partial_{x_1}a_2 V_7,V_1)_{L^2}\\
&+(\partial_{x_1}a_1 V_2,V_4)_{L^2}+(\partial_{x_1}a_2 V_7,V_4)_{L^2}\\
&+(\partial_{x_2}a_1 V_2,V_5)_{L^2}+(\partial_{x_2}a_2 V_7,V_5)_{L^2}\\
&+(\partial_{x_2}a_1 V_2,V_8)_{L^2}+(\partial_{x_2}a_2 V_7,V_8)_{L^2}.
\end{split}
\]
Making use of Glaeser's inequality we therefore have that $(S_1V,V)_{L^2}$ can be estimated by the sum of $\Vert V_1\Vert^2_{L^2}$,  $(a_1V_2,V_2)_{L^2}$, $(a_2V_7, V_7)_{L^2}$, $\Vert V_4\Vert_{L^2}^2$, $\Vert V_5\Vert^2_{L^2}$, $\Vert V_8\Vert^2_{L^2}$. Hence, we conclude that also $(S_1V,V)_{L^2}\le c'E(t)$ for some suitable constant $c'>0$. It follows that a combination of (H) and (LC) leads to 
\beq
\label{est_2_V}
|(({\wt{Q}}{\wt{B}}+{\wt{B}}^\ast{\wt{Q}})V,V)_{L^2}|\le c_2 E(t),
\eeq
for some $c_2>0$. Finally, to estimate $({\wt{Q}}V,{\wt{F}})_{L^2}$ we write it as 
\[
({\wt{Q}}V,\nabla_x {F})_{L^2} +({\wt{Q}}V, \left(
	\begin{array}{c}
	(\partial_{x_1}{B})U\\ 
	(\partial_{x_2}{B})U\\ 
	\vdots \\
	(\partial_{x_k}{B})U\\
	\vdots\\
	(\partial_{x_n}{B})U\\ 
	\end{array}
	\right))_{L^2}.
\]
It is immediate to see that 
\beq
\label{est_3_V}
2|({\wt{Q}}V,\nabla_x {F})_{L^2}|\le E(t)+\Vert f\Vert_{H^1}^2.
\eeq
To estimate the remaining term, that for brevity we will call $(\wt{Q}V, T_2)_{L^2}$, we argue as in \cite{G20} Proposition 4.9 (iii). We begin by noting that
\[
\partial_t U_j=V_{(j-1)(n+2)}.
\]
for $j\neq 1, n+2$. So, by the fundamental theorem of calculus combined with Cauchy-Schwarz and the Minkowski's inequality in integral form, we get
\[
\begin{split}
2|(\wt{Q}V, T_2)_{L^2}|&\le E(t)+ c(n, T, \max_{i=1,\dots,n, |\alpha|=1}(\Vert \partial^\alpha b_i\Vert_{\infty}^2, \Vert \partial^\alpha c\Vert_\infty^2, \Vert \partial^\alpha d\Vert_\infty^2))\\
&\biggl(\int_{0}^t \sum_{j=2}^{n+1}\Vert V_{(j-1)(n+2)(s)}\Vert^2_{L^2}\, ds+ \sum_{j=2}^{n+1}\Vert U_j(0)\Vert^2_{L^2}+\Vert U_1\Vert_{L^2}^2+\Vert U_{n+2}\Vert_{L^2}^2\biggr).
\end{split}
\]
Note that we already know how to estimate $\Vert U_1\Vert_{L^2}^2$ and $\Vert U_{n+2}\Vert_{L^2}^2$. It follows that if

\begin{itemize}
\item[(C)] \underline{the lower order terms have bounded first order derivatives},
\end{itemize}

then there exists a constant $c_3>0$ such that 
\[
2|(\wt{Q}V, T_2)_{L^2}|\le E(t)+c_3\biggl(\int_{0}^t E(s)\, ds+\Vert g_0\Vert_{H^1}^2+\Vert g_1\Vert_{L^2}^2+\int_{0}^t\Vert f(s)\Vert_{L^2}^2\, ds\biggr).
\]
Combining \eqref{est_3_V} with the estimate obtained above we conclude that under condition (C) there exists a constant $c_3>0$ such that 
\beq
\label{est_4_V}
|2({\wt{Q}}V,{\wt{F}})_{L^2}|\le E(t)+\Vert f\Vert_{H^1}^2+c_3\biggl(\int_{0}^t E(s)\, ds+\Vert g_0\Vert_{H^1}^2+\Vert g_1\Vert_{L^2}^2+\int_{0}^t\Vert f(s)\Vert_{L^2}^2\biggr).
\end{equation}
By collecting \eqref{est_1_V}, \eqref{est_2_V}, and \eqref{est_4_V} we have that under the hypotheses (H), (LC) and (C) the following estimate
\[
\begin{split}
 \frac{dE(t)}{dt} &\le (c_1+c_2+1) E(t)+\Vert f\Vert_{H^1}^2+c_3\biggl(\int_{0}^t E(s)\, ds+\Vert g_0\Vert_{H^1}^2+\Vert g_1\Vert_{L^2}^2+\int_{0}^t\Vert f(s)\Vert_{L^2}^2\biggr)\\
 &\le c'\biggl(E(t)+\int_{0}^t E(s)\, ds+\Vert g_0\Vert_{H^1}^2+\Vert g_1\Vert_{L^2}^2+\Vert f(t)\Vert_{H^1}^2+\Vert f\Vert_{L^\infty\times L^2}^2\biggr),
 \end{split}
\]
holds, where the constant $c'$ depends on $M$, $T$, $n$ and the $L^\infty$-norms of the first derivatives of the lower order terms. Now by applying a Gr\"onwall's type lemma (Lemma 6.2 in \cite{ST} or Lemma 4.10 in \cite{G20}) there exists a constant $C'>0$ depending exponentially on $M,n,T$ and the $L^\infty$-norms of the first derivatives of the lower order terms such that 
\beq
\label{est_E_V}
E(t)\le C'\biggl(E(0)+\int_0^t \Vert f(s)\Vert_{H^1}^2\, ds+\Vert g_0\Vert_{H^1}^2+\Vert g_1\Vert_{L^2}^2+\Vert f\Vert_{L^\infty\times L^2}^2\biggr).
\eeq
By definition of the energy $E(t)$ we have that 
\beq
\label{est_fin_tV}
\sum_{j=2}^{n+1}\Vert V_{(j-1)(n+2)}\Vert_{L^2}^2\le C'\biggl(\int_0^t \Vert f(s)\Vert_{H^1}^2\, ds+\Vert g_0\Vert_{H^2}^2+\Vert g_1\Vert_{H^1}^2+\Vert f\Vert_{L^\infty\times L^2}^2\biggr),
\eeq
where $C''>0$ depends on $C'$ and the $L^\infty$-norm of the coefficients $a_i$, $i=1,\cdots,n$. By using the relation $\partial_t U_j=V_{(j-1)(n+2)}$, for $j\neq 1, n+2$ and arguing as in Proposition 4.8 in \cite{G20} by fundamental theorem of calculus and Minkowski's integral inequality we can rewrite \eqref{est_fin_tV} in terms of the entries $U_j$. Hence, we conclude that for $j\neq 1, n+2$,
\[
\Vert U_j\Vert_{L^2}^2\le C\biggl(\int_0^t \Vert f(s)\Vert_{H^1}^2\, ds+\Vert g_0\Vert_{H^2}^2+\Vert g_1\Vert_{H^1}^2+\Vert f\Vert_{L^\infty\times L^2}^2\biggr),
\]
for a suitable constant $C>0$. As explained in \cite{G20} this leads to the $L^2$ well-posedness of the Cauchy problem in $U$.\\[0.2cm]
{\bf Summary.} We have proven that under the hypotheses
\begin{itemize}
\item[(H)] the coefficients $a_i$ are non-negative and bounded for all $i=1,\dots, n$ with bounded second order derivatives,
\item[(LC)] $b_i^2\prec a_i$ for all $i=1,\cdots,n$ and the lower order terms $c$ and $d$ are bounded,  
\item[(C)] the lower order terms $b_i, c, d$, $i=1,\cdots, n$ have bounded first order derivatives,
\end{itemize}
the Cauchy problem for the homogeneous wave equation
\[
\begin{split}
\partial_t^2u-\sum_{i=1}^na_i(x)\partial^2_{x_i} u+\sum_{i=1}^nb_i(x)\partial_{x_i}u+c(x)\partial_t u+d(x)u&=f(t,x),\\
u(0, x)&=g_0(x),\\
\partial_t u(0,x)&=g_1(x)\\
\end{split}
\]
on $[0,T]\times\R^n$ is well-posed with loss of derivatives, i.e., 
\[
\begin{split}
\Vert u(t)\Vert_{L^2}^2&\le  c''\biggl(\Vert g_0\Vert_{H^1}^2+\Vert g_1\Vert_{L^2}^2+\int_{0}^t\Vert f(s)\Vert_{L^2}^2\, ds\biggr),\\
\Vert \partial_t u(t)\Vert_{L^2}^2&\le  c''\biggl(\Vert g_0\Vert_{H^1}^2+\Vert g_1\Vert_{L^2}^2+\int_{0}^t\Vert f(s)\Vert_{L^2}^2\, ds\biggr),\\
\Vert u(t)\Vert_{H^1}^2&\le  c''\biggl(\Vert g_0\Vert_{H^2}^2+\Vert g_1\Vert_{H^1}^2+\int_{0}^t\Vert f(s)\Vert_{H^1}^2\, ds+\Vert f\Vert_{L^\infty\times L^2}^2\biggr).
\end{split}
\]

\begin{remark}
Note that one can iterate the previous scheme and get Sobolev estimates of every order for the solution $u$. However, this requires bigger and bigger matrices and the boundedness of more and more derivatives of the equation coefficients. The system that we get at the step $k$ has a block diagonal matrix generated by $A$ in the principal part, a matrix of lower order terms defined via $B$ and the first order $x$-derivatives of $A$ and the right-hand side which depends on the solutions of the previous $k-1$ steps, and $x$-derivatives of $B$ $A$ and $F$ up to order $k-1$. Therefore, by assuming that the equation coefficients belong to $B^\infty(\R^n)$ we can formulate the hypotheses (H) and (L) as   
\begin{itemize}
\item[(H)] the coefficients $a_i$ are non-negative for all $i=1,\dots, n$,
\end{itemize}
and the Levi conditions
\begin{itemize}
\item[(LC)] $b_i^2\prec a_i$ for all $i=1,\cdots,n$,
\end{itemize}
and our Cauchy problem is well-posed in every Sobolev space:
\[
\Vert u(t)\Vert_{H^k}^2\le  c_k\biggl(\Vert g_0\Vert_{H^{k+1}}^2+\Vert g_1\Vert_{H^k}^2+\int_{0}^t\Vert f(s)\Vert_{H^k}^2\, ds+\Vert f\Vert_{L^\infty\times H^{k-1}}^2\biggr),
\]
for all $k\in\N$. Note that when we deal with a homogeneous equation then by \cite{G20} and we get estimates as above without $\Vert f\Vert_{L^\infty\times H^{k-1}}^2$. This is consistent with the well-posedness result obtained in \cite{ST} for homogeneous hyperbolic equations with space-dependent coefficients in dimension 1.
 \end{remark}

\begin{remark}
The following table summarises the main steps employed above in obtaining Sobolev estimates of any order for the solution $u$. For the sake of simplicity, we assume $n=1$, but the same argument can be applied to any space dimension.

\begin{center}
\begin{tabular}{|c|c|c|c|}
\hline
System in & Hypotheses & $L^2$ estimates for & $u$ belongs to \\
\hline
$U$ & $a\ge 0$ and $(LC)$ & $U_1$ and $U_3$ & $L^2$\\
\hline
$V$ & $a\ge 0$, $(LC)$ and $B^\infty$ coefficients & $V_3$ and $U_2=V_1$ & $H^1$\\
\hline
$W$ & $a\ge 0$, $(LC)$ and $B^\infty$ coefficients & $W_3$ and $V_2=W_1$ & $H^2$\\
\hline
$\cdots$ & $\cdots$ & $\cdots$ & $\cdots$\\
\hline
\end{tabular}
\end{center}
\end{remark}
{\bf Conclusion}. We have proven, by symmetrisation method, the following $C^\infty$ well-posedness result.
\begin{theorem}
\label{theo_m_2}
Let
\[
\begin{split}
\partial_t^2u-\sum_{i=1}^na_i(x)\partial^2_{x_i} u+\sum_{i=1}^nb_i(x)\partial_{x_i}u+c(x)\partial_t u+d(x)u&=f(t,x),\\
u(0,x)&=g_0\in C^\infty_c(\R^n),\\
\partial_tu(0,x)&=g_1\in C^\infty_c(\R^n),
\end{split}
\]
where, $t\in[0,T]$, $x\in\R^n$, all the equation coefficients are real-valued and belong to $B^\infty(\R^n)$, $a_i\ge 0$ for all $i=1,\cdots, n$ and $f\in C([0,T], C^\infty_c(\R^n))$. Then, under the Levi conditions
\[
b_i^2\prec a_i, \qquad i=1,\cdots, n,
\]
the Cauchy problem is $C^\infty$ well-posed, i.e., there exists a unique solution 
\[
u\in C^2([0,T], C^\infty(\R^n)).
\]
\end{theorem}
We now ask ourselves if a similar result still holds for order $m>2$ and which Levi conditions to formulate on the lower order terms to guarantee well-posedness in every Sobolev space and therefore $C^\infty$ well-posedness.

\section{Third order hyperbolic equations in space dimension 1}
We begin by studying a third order hyperbolic equation in space dimension 1, i.e., we assume $x\in\R$ and $t\in[0,T]$. We want to investigate the well-posedness of the Cauchy problem 
\beq
\label{CP_3_toy_1}
\begin{split}
\partial_t^3u-a(x)\partial_t\partial^2_{x} u+b_1(x)\partial_{x}^2u+b_2(x)\partial_t\partial_{x}u+b_3(x)\partial^2_tu&=f(t,x),\\ 
u(0,x)&=g_0(x),\\
\partial_tu(0,x)&=g_1(x),\\
\partial_t^2u(0,x)&=g_2(x),
\end{split}
\eeq
where 
$a(x)\ge 0$ for all $x\in\R$ and all the equation coefficients are real-valued and belongs to $B^\infty(\R)$. We will also assume that the initial data are smooth functions with compact support and that $f\in C([0,T], C^\infty_c(\R))$.

As for second order equations we will employ symmetrisation techniques. For the sake of simplicity we work with lower order terms of order 2 to allow a system transformation which can be easily adapted to higher space dimensions. Setting
\[
U=(\partial_{x}^2 u,\partial_{x}\partial_t u,\partial_t^2u)^T.
\]
we can rewrite \eqref{CP_3_toy_1} as 
\beq
\label{CP_3_U_1}
\begin{split}
\partial_t U&=A(x)\partial_xU+B(x)U+F,\\
U(0,x)&=U_0(x)=(g^{(2)}_0,g^{(1)}_1,g_2)^T.
\end{split}
\eeq
where the matrices $A$ and $B$ have size $3\times 3$ and $F$ is the -column $(0,0,f)^T$.
The matrix $A$ is in Sylvester form and the matrix $B$ of the lower order terms has only the last row non identically zero. In detail,
\[
A=\left(
	\begin{array}{ccc}
	0& 1 & 0\\   
	0& 0 & 1\\  
	0 & a & 0  
	\end{array}
	\right),
B=\left(
	\begin{array}{ccc}
	0 & 0 & 0\\ 
	0 & 0 & 0\\
	-b_1 & -b_2 & -b_3 
  \end{array}
	\right).
\]
From the general theory of Sylvester matrices and symmetrisation (see the Appendix) we have that the matrix  
	\[
	Q=\frac{1}{3}\left(
	\begin{array}{ccc}
	a^2 & 0 & -a\\
	0 & 2a & 0\\
	-a & 0 & 3
	\end{array}
	\right).
	\]
is the standard symmetriser of $A$. Indeed,
\[
QA= A^\ast Q=  \frac{2}{3}\left(
	\begin{array}{ccc}
	0 & 0 & 0\\
	0 & 0 & a\\
	0 & a & 0
	\end{array}
	\right).
\]
If we denote the roots of our equation with $\lambda_1=-\sqrt{a}$, $\lambda_2=0$ and $\lambda_3=\sqrt{a}$ we have that 
\[
{\rm det}Q=\frac{1}{27}\sum_{1\le i<j\le 3}(\lambda_i-\lambda_j)^2=\frac{8a^3}{27}.
\]
 Note that, differently from the case $m=2$, the symmetriser $Q$ is not diagonal but it is nearly diagonal.  This means that, given the diagonal matrix
\[
\Psi=\left(
	\begin{array}{ccc}
	a^2 & 0 & 0\\
	0 & 2a & 0\\
	0 & 0 & 1
	\end{array}
	\right)
\]
we can find suitable constants $c_1,c_2>0$ such that 
\[
c_1\lara{\Psi v,v}\le \lara{Qv,v}\le c_2\lara{\Psi v,v}
\]
for all $v\in\R^3$. Since the matrices involved depend on $x$, the inequality holds uniformly in $x\in\R$. Indeed,
\[
\begin{split}
3\lara{Qv,v}&= \Vert av_1\Vert^2+2\lara{av_2,v_2}+3\Vert v_3\Vert^2-2\lara{av_1,v_3}\\
&\le  \Vert av_1\Vert^2+2\lara{av_2,v_2}+3\Vert v_3\Vert^2+\Vert av_1\Vert^2+\Vert v_3\Vert^2\\
&\le 4\lara{\Psi v,v}\\
3\lara{Q v,v}&= \Vert av_1\Vert^2-2\lara{av_1,v_3}+\Vert v_3\Vert^2+2\lara{av_2,v_2}+2\Vert v_3\Vert^2\\
&=\Vert av_1-v_3\Vert^2+\Vert v_3\Vert^2+2\lara{av_2,v_2}+\Vert v_3\Vert^2\\
&\ge \frac{1}{2}\Vert av_1\Vert^2+2\lara{av_2,v_2}+\Vert v_3\Vert^2\\
&\ge \frac{1}{2}\lara{\Psi v,v}.
\end{split}
\]
It follows that $c_1=\frac{1}{6}$ and $c_2=\frac{4}{3}$. We can now define the energy
\[
E(t)=(Q(t)U,U)_{L^2},
\]
and proceed with the energy estimates like for the wave equation in the previous section and in \cite{G20}.  We end up with  
 \[
 \begin{split}
 &\frac{dE(t)}{dt}=(\partial_t(QU),U)_{L^2}+(QU,\partial_tU)_{L^2}\\
 &=-((QA)'U,U)_{L^2}+((QB+B^\ast Q)U,U)_{L^2}+2(QU,F)_{L^2}.
\end{split}
 \]
 Let us focus on the second term: $((QB+B^\ast Q)U,U)_{L^2}$. We want to find suitable Levi conditions on $B$ such that we can bound  $((QB+B^\ast Q)U,U)_{L^2}$ with the energy $E(t)$. This will lead to the $L^2$ well-posedness of the corresponding Cauchy problem and by iteration to the $C^\infty$ well-posedness. In other words we want to guarantee that 
 \[
  ((QB+B^\ast Q)U,U)_{L^2}\le c\, E(t),
 \]
 for some constant $c>0$ uniformly in $t\in[0,T]$ and $x\in\R$.  Making use of the fact that $Q$ is nearly diagonal, it suffices to prove that there exists $c>0$ such that 
 \[
  ((QB+B^\ast Q)U,U)_{L^2}\le c(\Psi U,U)_{L^2}.
 \]
 By straightforward computations we have 
 \[
 QB+B^\ast Q=\frac{1}{3}\left(
	\begin{array}{ccc}
	2ab_1 & ab_2 & ab_3-3b_1\\
	ab_2 & 0 & -3b_2\\
	ab_3-3b_1 & -3b_2 & -6b_3
	\end{array}
	\right)
 \]
 and therefore
 \[
 \begin{split}
 3 ((QB+B^\ast Q)U,U)_{L^2}&=(2ab_1U_1,U_1)_{L^2}+2((ab_3-3b_1)U_1,U_3)_{L^2}-6(b_2U_2,U_3)_{L^2}\\
& +(2ab_2U_1,U_2)_{L^2}-6(b_3U_3,U_3)_{L^2}.
 \end{split}
 \]
If we now compare this term with 
\[
(\Psi U,U)_{L^2}=(a^2U_1,U_1)_{L^2}+(2aU_2,U_2)_{L^2}+(U_3,U_3)_{L^2}
\]
we immediately see that if 
\beq
\label{LC_3}
\begin{split}
|b_1(x)|&\le c_1 a(x),\\
|b_2(x)|&\le c_2\sqrt{a(x)},\\
|b_3(x)|&\le c_3
\end{split}
\eeq
uniformly in $x\in\R$ then there exists a constant $c>0$ such that
\[
 3 ((QB+B^\ast Q)U,U)_{L^2}\le c(\Psi U,U)_{L^2}\le 6c E(t),
\]
as desired. Clearly, since the equation coefficients are bounded the third Levi condition is redundant. 

The Levi conditions above are sufficient to ensure $L^2$ well-posedness. Inspired by \cite{G20}, we now derive the system \eqref{CP_3_toy_1} with respect to $x$ and proceed with the energy estimates. This will lead to an extra Levi condition which will guarantee $H^1$ well-posedness. In detail, for $V=\partial_xU$ we get that if $U$ solves \eqref{CP_3_U_1} then $V$ solves
\beq
\label{CP_3_V_1}
\begin{split}
\partial_t V&=A(x)\partial_xV+(A'(x)+B(x))V+B'(x)U+\partial_xF,\\
V(0,x)&=(g^{(3)}_0,g^{(2)}_1,g^{(1)}_2)^T.
\end{split}
\eeq
The energy for this system is still the same since the principal part has not changed. However, we have a new matrix of lower order terms and a new right-hand side. It follows that we need to estimate 
\[
((QA'+(A')^\ast Q)V,V)_{L^2}+((QB+B^\ast Q)V,V)_{L^2}
\]
with 
\[
E(t)=(QV,V)_{L^2}.
\]
The Levi conditions \eqref{LC_3} guarantee that $((QB+B^\ast Q)U,U)_{L^2}$ is bounded by the energy. By direct computations, we immediately see that 
\[
3((QA'+(A')^\ast Q)V,V)_{L^2}=-2(aa'V_1,V_2)_{L^2}+2(3a'V_2,V_3).
\]
Hence
\[
3|(QA'+(A')^\ast Q)V,V)_{L^2}|\le (a^2V_1,V_1)_{L^2}+\Vert a'V_2\Vert_{L^2}^2+\Vert a'V_2\Vert_{L^2}^2+\Vert 3V_3\Vert_{L^2}^2,
\]
and, since by the Glaeser's inequality $|a'(x)|^2\le 2Ma(x)$, we have that
\beq
\label{est_import_3}
\begin{split}
3|(QA'+(A')^\ast Q)V,V)_{L^2}|&\le (a^2V_1,V_1)_{L^2}+2M\Vert aV_2\Vert_{L^2}^2+2M\Vert aV_2\Vert_{L^2}^2+\Vert 3V_3\Vert_{L^2}^2\\
&\le c(M)(\Psi V,V)_{L^2}\\
&\le 6c(M)E(t).
\end{split}
\eeq
It remains to estimate the term $B'U$ with the energy $E(t)$. Note that 
\[
B'U=\left(
	\begin{array}{c}
	0 \\
	0 \\
	 -b_1'U_1-b_2'U_2-b_3'U_3
        \end{array}
	\right),
\]
and, from the previous analysis we can estimate $\Vert U_3\Vert_{L^2}^2$ in terms of $f$ and the initial data of \eqref{CP_3_toy_1}. The only terms that we need to estimate are therefore $\Vert -b_1'U_1\Vert^2_{L^2}$ and $\Vert -b_2'U_2\Vert^2_{L^2}$. We begin by noting that $\partial_tU_1=V_2$ and $\partial_t U_2=V_3$. Arguing as in \cite{G20} (Subsection 3.3) we can write
\begin{multline}
 \label{est_E_U_1}
\Vert -b'_1U_1\Vert^2_{L^2}=\biggl\Vert\int_0^t b_1'\partial_tU_1(s)\, ds+b_1'U_1(0)\biggr\Vert^2\le 2\biggl\Vert \int_0^t b_1'\partial_tU_1(s)\, ds\biggr\Vert^2_{L^2}+2\Vert b_1'U_1(0)\Vert_{L^2}^2\\
\le2\biggl(\int_0^t\Vert b_1'V_2(s)\Vert_{L^2}\, ds\biggr)^2+2\Vert b_1'U_1(0)\Vert_{L^2}^2
\le 2T\int_0^t((b_1')^2V_2,V_2)_{L^2}\, ds+2\Vert  b_1'U_1(0)\Vert_{L^2}^2\\
\le 2T\int_0^t((b_1')^2V_2,V_2)_{L^2}\, ds+2\Vert b_1'\Vert_{L^\infty}^2\Vert g_0\Vert_{H^2}^2
\end{multline}
and
\begin{multline}
 \label{est_E_U_2}
\Vert -b_2'U_2\Vert^2_{L^2}=\biggl\Vert\int_0^t b_2'\partial_tU_2(s)\, ds+b_2'U_2(0)\biggr\Vert^2\le 2\biggl\Vert \int_0^t b_2'\partial_tU_2(s)\, ds\biggr\Vert^2_{L^2}+2\Vert b_2'U_2(0)\Vert_{L^2}^2\\
\le2\biggl(\int_0^t\Vert b_2'V_3(s)\Vert_{L^2}\, ds\biggr)^2+2\Vert b_2'U_2(0)\Vert_{L^2}^2
\le 2T\int_0^t((b_2')^2V_3,V_3)_{L^2}\, ds+2\Vert  b_2'U_2(0)\Vert_{L^2}^2\\
\le 2T\int_0^t((b_2')^2V_3,V_3)_{L^2}\, ds+2\Vert b_2'\Vert_{L^\infty}^2\Vert g_1\Vert_{H^1}^2.
\end{multline}
Combining \eqref{est_E_U_1} with \eqref{est_E_U_2} it is clear that if there exists $c_4>0$ such that 
\beq
\label{LC_4}
|b_1'(x)|\le c_4\sqrt{a(x)},
\eeq
for all $x\in\R$ then
\[
\Vert -b_1'U_1\Vert^2_{L^2}+\Vert -b_2'U_2\Vert^2_{L^2}\le c\biggl(\int_0^t E(s)\, ds+\Vert g_0\Vert_{H^2}^2+\Vert g_1\Vert_{H^1}^2\biggr),
\]
where the constant $c>0$ depends on $c_4$, $T$ and the $L^\infty$-norms of $b_1'$ and $b_2'$. We therefore conclude, by applying a Gr\"onwall type lemma as in \cite[Lemma 4.10]{G20},  that in order to get $H^1$-estimates for our solution $U$ we need to add a Levi condition on the lower order terms, namely on the first order derivative of $b_1$. Summarising,
\[
\begin{split}
|b_1(x)|&\le c_1 a(x),\\
|b_2(x)|&\le c_2\sqrt{a(x)},\\
|b_3(x)|&\le c_3,\\
|b_1'(x)|&\le c_4\sqrt{a(x)}.
\end{split}
\]
Note that these Levi conditions are enough to guarantee estimates in the next Sobolev order. Indeed, by taking an extra $x$ derivative we get, for $W=\partial_xV$, the Cauchy problem
\beq
\label{CP_3_W_1}
\begin{split}
\partial_t W&=A(x)\partial_xW+(2A'(x)+B(x))W+(A''(x)+2B'(x)V+B''(x)U+\partial^2_xF,\\
W(0,x)&=(g^{(4)}_0,g^{(3)}_1,g^{(2)}_2)^T.
\end{split}
\eeq
Glaeser's inequality and the first three Levi conditions allow us to estimate $(2A'(x)+B(x))W$ with the energy $E(t)$. Arguing as in \eqref{est_E_U_1} and \eqref{est_E_U_2} (replace $U$ and $V$ with $V$ and $W$, respectively) and making use of the fact that $\partial_tV_1=W_2$ and $\partial_t V_2=W_3$ we have that $\Vert A''(x)V\Vert_{L^2}^2$ can be estimated by $\int_0^t E(s)\, ds$ and the norms of the initial data. Since
\[
\Vert B'(x)V\Vert_{L^2}^2=\Vert b_1'V_1\Vert_{L^2}^2+\Vert b_2'V_2\Vert_{L^2}^2+\Vert b_3'V_3\Vert_{L^2}^2,
\]
we easily see that the fourth Levi condition on $b'_1$ is needed to estimate $\Vert b_1'V_1\Vert_{L^2}^2$ and therefore $\Vert B'(x)V\Vert_{L^2}^2$ in terms of $\int_0^t E(s)\, ds$ and suitable Sobolev norms of the initial data. We are now ready to prove the following theorem. Note that since the coefficients are bounded the Levi condition on $b_3$ is automatically fulfilled.
\begin{theorem}
\label{theo_3_case_1}
Let
\beq
\label{CP_3_no_l_1}
\begin{split}
\partial_t^3u-a(x)\partial_t\partial^2_{x} u+b_1(x)\partial_{x}^2u+b_2(x)\partial_t\partial_{x}u+b_3(x)\partial^2_tu&=f(t,x),\\ 
u(0,x)&=g_0(x),\\
\partial_tu(0,x)&=g_1(x),\\
\partial_t^2u(0,x)&=g_2(x),
\end{split}
\eeq
where $a\ge 0$ and all the equation coefficients and are real-valued and belong to $B^\infty(\R)$. Let $f\in C^3([0,T], H^\infty(\R))$.
\begin{itemize}
\item[(i)] Under the Levi conditions $(LC)$:
\[
|b_1|\prec a,\quad |b_2|\prec \sqrt{a},
\]
the Cauchy problem has a unique solution in $C^3([0,T], L^2(\R))$ provided that $g_0\in H^2(\R)$, $g_1\in H^1(\R)$ and $g_2\in L^2(\R)$. Moreover,
\[
\Vert u(t)\Vert_{L^2}^2\le c\biggl(\int_{0}^t\Vert f(s)\Vert_{L^2}^2\, ds+\Vert g_0\Vert_{H^2}^2+\Vert g_1\Vert_{H^1}^2+\Vert g_2\Vert_{L^2}^2 \biggr).
\]
\item[(ii)] Under the Levi conditions $(LC)_1$:
\[
|b_1|\prec a,\quad |b_2|\prec \sqrt{a},\quad |b'_1|\prec \sqrt{a},
\]
the Cauchy problem has a unique solution in $C^3([0,T], H^1(\R))$ provided $g_0\in H^3(\R)$, $g_1\in H^2(\R)$ and $g_2\in H^1(\R)$. Moreover,
\[
\Vert u\Vert_{H^1}^2\le c\biggl(\int_0^t \Vert f(s)\Vert_{H^1}^2\, ds+\Vert f\Vert_{L^\infty\times L^2}^2+\Vert g_0\Vert _{H^3}^2+\Vert g_1\Vert_{H^2}^2+\Vert g_2\Vert_{H^1}^2 \biggr).
\]
\item[(iii)] Finally under the Levi conditions $(LC)_1$, the Cauchy problem has a unique solution in $C^3([0,T], H^k(\R))$, $k\in\N$, provided $g_0\in H^{k+2}(\R)$, $g_1\in H^{k+1}(\R)$ and $g_2\in H^k(\R)$. Moreover,
\[
\Vert u\Vert_{H^k}^2\le c\biggl(\int_0^t \Vert f(s)\Vert_{H^k}^2\, ds+\Vert f\Vert_{L^\infty\times H^{k-1}}^2+\Vert g_0\Vert _{H^{k+2}}^2+\Vert g_1\Vert_{H^{k+1}}^2+\Vert g_2\Vert_{H^k}^2 \biggr).
\]

\end{itemize}
\end{theorem}
It follows immediately that our Cauchy problem is $C^\infty$ well-posed.
\begin{theorem}
\label{theo_m_3_1}
Let $f\in C^3([0,T], C^\infty_c(\R))$ and $g_0,g_1,g_2\in C^\infty_c(\R)$. Then, the Cauchy problem \eqref{CP_3_no_l_1} is $C^\infty$ well-posed provided that $a\ge 0$, all the equation coefficients are real-valued and belong to $B^\infty(\R)$ and the Levi conditions 
\[
|b_1|\prec a,\quad |b_2|\prec \sqrt{a},\quad |b'_1|\prec \sqrt{a}
\]
are fulfilled.
\end{theorem}
\begin{remark}
Comparing our result with the one obtained in \cite{ST21} we see that our Levi conditions are more general in the sense that they replace equalities with bounds. Indeed, according to \cite{ST21}, $C^\infty$ well-posedness is obtained when the polynomial of the lower order terms has a proper decomposition with respect to the principal part. This means that 
\[
b_1(x)+b_2(x)\tau+b_3(x)\tau^2=\sum_{k=1}^3 l_k(x)P_{\widehat{k}}(x,\tau), 
\]
where $l_k$ are bounded functions and 
\[
P_{\widehat{k}}(x,\tau)=\Pi_{j=1,2,3, j\neq k}(\tau-\tau_k(x)), \qquad \tau_1=0,\, \tau_2=-\sqrt{a(x)},\, \tau_3=+\sqrt{a(x)}.
\]
By direct computations, it follows that the coefficients $b_i$ are determined by the principal part $a$ as  $b_1(x)=\lambda_1(x)a(x)$ and $b_2(x)=\lambda_2(x)\sqrt{a(x)}$, with $\lambda_i$ bounded, $i=1,2$. Working under the assumption that the equation coefficients are elements of $B^\infty(\R)$ we can assume that $\lambda_1'$ is bounded as well. It follows that the Levi conditions in \cite{ST21} imply the ones in Theorem \ref{theo_m_3_1}, since, by Glaeser's inequality, 
\[
|b_1'|=|\lambda_1'a+\lambda_1a'|\prec a+|a'|\prec \sqrt{a}.
\]
 \end{remark}

\begin{proof}[Proof of Theorem \ref{theo_3_case_1}]
The existence of the solution $u$ is guaranteed by a perturbation argument (Nuji's approximation) similar to the one employed for the wave equation. We focus here on estimating the solution $u$ in terms of the initial data and the right-hand side. It is not restrictive to assume that the initial data are compactly supported and, by finite speed propagation, to assume that the solution is compactly supported with respect to $x$. This assumption can be later removed by density argument. We wok on the system in $U$ defined in 
\eqref{CP_3_U_1} and on the energy $E(t)$ defined by the symmetriser which leads to
 \[
 \frac{dE(t)}{dt}=-((QA)'U,U)_{L^2}+((QB+B^\ast Q)U,U)_{L^2}+2(QU,F)_{L^2}.
 \]
By straightforward computations we have that 
\[
((QA)'U,U)_{L^2}=\frac{4}{3}(a'U_2,U_3)_{L^2}.
\]
Since by Glaeser's inequality $|a'(x)|^2\le 2M a(x)$, where $\Vert a''\Vert_{L^\infty}\le M$ we conclude that  
\[
|(a'U_2,U_3)_{L^2}|\le \Vert a'U_2\Vert^2_{L^2}+\Vert U_3\Vert^2_{L^2}\le 3M(\frac{2}{3}a_2U_2,U_2)_{L^2}+\Vert U_3\Vert^2_{L^2}\le \max\{3M,1\}E(t).
\]
As observed in the arguments leading to \eqref{LC_3} under the Levi conditions (LC) the term $(QB+B^\ast Q)U,U)_{L^2}$ is bounded by the energy $E(t)$. Finally, 
\[
|3(QU, F)_{L^2}|=|(-aU_1+3U_3,f)_{L^2}|\le (a^2U_1,U_1)_{L^2}+\Vert f\Vert_{L^2}^2+9\Vert U_3\Vert^2_{L^2}+\Vert f\Vert_{L^2}^2.
\]
By the fact that $Q$ is nearly diagonal, i.e,
\[
(QU,U)_{L^2}\ge \frac{1}{6}(\Psi U,U)_{L^2},
\]
we easily conclude that $(QU, F)_{L^2}$ can be bounded by $E(t)+\Vert f\Vert^2_{L^2}$. Hence, there exist a constant $C>0$ depending on the $L^\infty$-norm of the second derivative of $a$ and the Levi conditions $(LC)$ such that 
\[
 \frac{dE(t)}{dt}\le C(E(t)+\Vert f(t)\Vert_{L^2}^2).
\]
By Gr\"onwall's lemma and the bound from below 
\[
\Vert \partial_t u\Vert_{L^2}^2=\Vert U_3\Vert_{L^2}^2\le E(t)
\] 
we obtain the following estimate for a suitable constant $c>0$:
\[
\begin{split}
\Vert \partial_tu(t)\Vert_{L^2}^2\le E(t)&\le c\biggl(E(0)+\int_{0}^t\Vert f(s)\Vert_{L^2}^2\, ds\biggr)\\
&\le  c\biggl(\Vert g_0\Vert_{H^2}^2+\Vert g_1\Vert_{H^1}^2+\Vert g_2\Vert_{L^2}^2+\int_{0}^t\Vert f(s)\Vert_{L^2}^2\, ds\biggr).
\end{split}
\]
Note that here we have also used the fact that the coefficients $a$ is bounded, with bounded derivatives of any order. 

To estimate the other entries of $U$ we pass to the system \eqref{CP_3_V_1} in $V$ obtained by differentiating once more with respect to $x$. In detail,
\[
\begin{split}
\partial_t V&=A(x)\partial_xV+(A'(x)+B(x))V+B'(x)U+\partial_xF,\\
V(0,x)&=(g^{(3)}_0,g^{(2)}_1,g^{(1)}_2)^T
\end{split}
\]
and 
 \[
 \begin{split}
\frac{dE(t)}{dt}
 &=-((QA)'V,V)_{L^2}+((QA'+(A')^\ast Q)V,V)_{L^2}+((QB+B^\ast Q)V,V)_{L^2}\\
 &+2(QV, B'U+\partial_xF)_{L^2}.
 \end{split}
 \]
As observed in the section leading to the formulation of the Levi conditions $(LC)$ we can estimate the first three addenda on the right-hand side of the formula above using the hypothesis on $a$ and $(LC)$. The extra Levi condition on $b_1'$ is needed to estimate $(QV, B'U)_{L^2}$. We get 
\beq
\label{term_import_3}
(QV, B'U)_{L^2}=\frac{1}{3}(-aV_1+3V_3, -b_1'U_1-b_2'U_2-b_3'U_3)_{L^2}
\eeq
where 
\begin{multline*}
|(-aV_1+3V_3, -b_1'U_1-b_2'U_2-b_3'U_3)_{L^2}| \prec\\
((a^2V_1,V_1)_{L^2}+\Vert V_3\Vert_{L^2}^2+\Vert b_1'U_1\Vert^2_{L^2}+\Vert b_2'U_2\Vert^2_{L^2}+\Vert b_3'U_3\Vert^2_{L^2}).
\end{multline*}
Note that 
\[
(a^2V_1,V_1)_{L^2}+\Vert V_3\Vert_{L^2}^2\le E(t)
\]
and, by the boundedness of $b_3'$, 
\[
\Vert b_3'U_3\Vert^2_{L^2}\prec \Vert U_3\Vert_{L^2}^2\prec\biggl(\Vert g_0\Vert_{H^2}^2+\Vert g_1\Vert_{H^1}^2+\Vert g_2\Vert_{L^2}^2+\int_{0}^t\Vert f(s)\Vert_{L^2}^2\, ds\biggr).
\]
Under the Levi condition $(LC)_1$ we obtain \eqref{est_E_U_1} which combined with \eqref{est_E_U_2} yields
\[
\Vert -b_1'U_1\Vert^2_{L^2}+\Vert -b_2'U_2\Vert^2_{L^2}\prec \biggl(\int_0^t E(s)\, ds+\Vert g_0\Vert_{H^2}^2+\Vert g_1\Vert_{H^1}^2\biggr).
\]
Concluding, there exists a constant $C>0$ such that 
\[
\begin{split}
\frac{dE(t)}{dt}&\le C\biggl(E(t)+\Vert g_0\Vert_{H^2}^2+\Vert g_1\Vert_{H^1}^2+\Vert g_2\Vert_{L^2}^2+\int_{0}^t\Vert f(s)\Vert_{L^2}^2\, ds+\int_0^t E(s)\, ds+\Vert f(t)\Vert_{H^1}^2\biggr)\\
&\le C\biggl(E(t)+\int_0^t E(s)\, ds+\Vert f(t)\Vert_{H^1}^2+\Vert f\Vert_{L^\infty\times L^2}^2+\Vert g_0\Vert _{H^2}^2+\Vert g_1\Vert_{H^1}^2+\Vert g_2\Vert_{L^2}^2 \biggr).
\end{split}
\]
Note that by comparison of $Q$ with $\Psi$ and since the coefficient $a$ is bounded we get that 
\[
E(0)\prec \Vert g_0\Vert _{H^3}^2+\Vert g_1\Vert_{H^2}^2+\Vert g_2\Vert_{H^1}^2.
\]
Hence, by application of a Gr\"onwall type Lemma (Lemma 4.10 in \cite{G20}), we obtain
\[
\Vert V_3\Vert_{L^2}^2\le C'\biggl(\int_0^t \Vert f(s)\Vert_{H^1}^2\, ds+\Vert f\Vert_{L^\infty\times L^2}^2+\Vert g_0\Vert _{H^3}^2+\Vert g_1\Vert_{H^2}^2+\Vert g_2\Vert_{H^1}^2 \biggr).
\]
Noting that $V_3=\partial_t U_2$ arguing as for the wave equation we get a similar estimate for $U_2$, i.e., 
\[
\Vert U_2\Vert_{L^2}^2\le C'\biggl(\int_0^t \Vert f(s)\Vert_{H^1}^2\, ds+\Vert f\Vert_{L^\infty\times L^2}^2+\Vert g_0\Vert _{H^3}^2+\Vert g_1\Vert_{H^2}^2+\Vert g_2\Vert_{H^1}^2 \biggr).
\]
Summarising, so far we have proven that under the Levi conditions $(LC)$
\beq
\label{est_u_3_0}
\Vert u(t)\Vert_{L^2}^2\le c\biggl(\Vert g_0\Vert_{H^2}^2+\Vert g_1\Vert_{H^1}^2+\Vert g_2\Vert_{L^2}^2+\int_{0}^t\Vert f(s)\Vert_{L^2}^2\, ds\biggr),
\eeq
and under the Levi conditions $(LC)_1$,
\beq
\label{est_u_3_1}
\Vert u\Vert_{H^1}^2\le c\biggl(\int_0^t \Vert f(s)\Vert_{H^1}^2\, ds+\Vert f\Vert_{L^\infty\times L^2}^2+\Vert g_0\Vert _{H^3}^2+\Vert g_1\Vert_{H^2}^2+\Vert g_2\Vert_{H^1}^2 \biggr).
\eeq
To estimate $U_1$ we need to derive once more the system in $V$ and use $W=\partial_x V$. By construction $\partial_t V_2=W_3$ and $V_2=\partial_t U_1$. So from $W_3$ we obtain an estimate for $V_2$ which can then be transferred to $U_1$. In detail,
\[
\begin{split}
\partial_t W&=A(x)\partial_xW+(2A'(x)+B(x))W+(2B'+A'')V+B''(x)U+\partial^2_xF,\\
W(0,x)&=(g^{(4)}_0,g^{(3)}_1,g^{(2)}_2)^T.
\end{split}
\]
This system has a very similar structure to the one in $V$ so our hypotheses will work perfectly on the principal part of the system and the matrix of order $0$. We want to give a closer look to the term 
\[
(QW, (2B'+A'')V+B''(x)U+\partial^2_xF)_{L^2}.
\]
By straightforward computations we get the quantity
\begin{multline*}
(-aW_1+3W_3, -2b_1'V_1-2b_2'V_2-2b_3'V_3+a''V_2)_{L^2}
+(-aW_1+3W_3,-b_1''U_1-b_2''U_2-b_3''U_3)_{L^2}\\
+(-aW_1+3W_3, \partial^2_x f)_{L^2}.
\end{multline*}
Recall that $\partial_t V_1=W_2$, $\partial_t V_2=W_3$. Hence, in analogy to \eqref{est_E_U_1} and \eqref{est_E_U_2}  using the fact that the coefficients are in $B^\infty(\R)$ and the Levi conditions $(LC)_1$ we get
\[
\begin{split}
(-aW_1+3W_3, -2b_1'V_1)_{L^2}&\prec E(t)+\int_0^t E(s)\, ds+\Vert W_2(0)\Vert^2_{L^2},\\
(-aW_1+3W_3, -2b_2'V_2+a''V_2)_{L^2}&\prec E(t)+\int_0^t E(s)\, ds+\Vert W_3(0)\Vert^2_{L^2},\\
(-aW_1+3W_3, -2b_3'V_3)_{L^2}&\prec E(t)+\Vert V_3\Vert_{L^2}^2\\
(-aW_1+3W_3,-b_2''U_2-b_3''U_3)_{L^2}&\prec E(t)+\Vert U_2\Vert_{L^2}^2+\Vert U_3\Vert_{L^2}^2,\\
(-aW_1+3W_3, \partial^2_x f)_{L^2}&\prec E(t)+\Vert f\Vert^2_{H^2}.
\end{split}
\]
Making use of the previous estimates on $U_2$, $U_3$ and $V_3$ and the initial data we obtain
\[
\begin{split}
(-aW_1+3W_3, -2b_1'V_1)_{L^2}&\prec E(t)+\int_0^t E(s)\, ds+\Vert g_1\Vert^2_{H^3},\\
(-aW_1+3W_3, -2b_2'V_2+a''V_2)_{L^2}&\prec E(t)+\int_0^t E(s)\, ds+\Vert g_2\Vert^2_{H^2},\\
(-aW_1+3W_3, -2b_3'V_3)_{L^2}&\prec E(t)+\int_0^t \Vert f(s)\Vert_{H^1}^2\, ds+\Vert f\Vert_{L^\infty\times L^2}^2\\
&+\Vert g_0\Vert _{H^3}^2+\Vert g_1\Vert_{H^2}^2+\Vert g_2\Vert_{H^1}^2\\
(-aW_1+3W_3,-b_2''U_2-b_3''U_3)_{L^2}&\prec E(t)+\int_0^t \Vert f(s)\Vert_{H^1}^2\, ds+\Vert f\Vert_{L^\infty\times L^2}^2\\
&+\Vert g_0\Vert _{H^3}^2+\Vert g_1\Vert_{H^2}^2+\Vert g_2\Vert_{H^1}^2,\\
(-aW_1+3W_3, \partial^2_x f)_{L^2}&\prec E(t)+\Vert f\Vert^2_{H^2}.
\end{split}
\]
It remains to estimate $(-aW_1+3W_3,-b_1''U_1)_{L^2}$. We make use of Minkowski's integral inequality and the relations $\partial_t U_1=V_2$, $\partial_t V_2=W_3$. We have 
\[
(-aW_1+3W_3,-b_1''U_1)_{L^2}\prec E(t)+\Vert b_1''U_1\Vert_{L^2}^2
\]
so we need to apply Minkowski's integral inequality to $\Vert b_1''U_1\Vert_{L^2}^2$ twice to be able to pass from $U_1$ to $W_3$. In detail, making also use of the fact that the equation coefficients have bounded derivatives of any order, we can write
\[
\begin{split}
\Vert b_1''U_1\Vert_{L^2}^2&=\biggl\Vert \int_0^t b''_1V_2(s)\, ds+ b_1''U_1(0)\biggr\Vert^2_{L^2}\\
&\prec \int_0^t \Vert b_1''V_2(s)\Vert^2_{L^2}\, ds+\Vert U_1(0)\Vert_{L^2}^2 \\
&=\int_0^t\biggl\Vert b_1''\biggl(\int_0^s W_3(r)\, dr+V_2(0)\biggr)\biggr\Vert_{L^2}^2 ds+\Vert U_1(0)\Vert_{L^2}^2 \\
&\prec \int_0^t\biggl\Vert b_1''\int_0^s W_3(r)\, dr\biggr\Vert_{L^2}^2 ds+\Vert V_2(0)\Vert_{L^2}^2 +\Vert U_1(0)\Vert_{L^2}^2\\
&\prec \int_0^t\int_{0}^s \Vert b_1'' W_3(r)\Vert_{L^2}^2\, dr\, ds+\Vert V_2(0)\Vert_{L^2}^2 +\Vert U_1(0)\Vert_{L^2}^2\\
&\prec \int_0^t \int_0^s E(r)\, dr\, ds+\Vert V_2(0)\Vert_{L^2}^2 +\Vert U_1(0)\Vert_{L^2}^2\\
&\prec \int_0^t E(s)\, ds+\Vert V_2(0)\Vert_{L^2}^2 +\Vert U_1(0)\Vert_{L^2}^2,
\end{split}
\]
where the constants hidden in the previous inequalities depend on the $L^\infty$-norms of the equation coefficients and the interval $[0,T]$. It follows that 
\[
(-aW_1+3W_3,-b_1''U_1)_{L^2}\prec E(t)+\int_0^t E(s)\, ds+\Vert g_1\Vert_{H^2}^2+\Vert g_0\Vert_{H^2}^2.
\]
Concluding, we have proven that 
\[
\begin{split}
&(QW, (2B'+A'')V+B''(x)U+\partial^2_xF)_{L^2}\prec E(t)+\int_0^t E(s)\, ds+\int_0^t \Vert f(s)\Vert_{H^1}^2\, ds\\
&+\Vert f\Vert_{H^2}^2+\Vert f\Vert^2_{L^\infty\times L^2}+\Vert g_0\Vert^2_{H^3}+\Vert g_1\Vert^2_{H^3}+\Vert g_2\Vert^2_{H^2}.
\end{split}
\]
We can now combine this estimate with the ones for the principal part of the system and the matrix of lower order terms and apply the Gr\"onwall type lemma. We get
\[
\begin{split}
\frac{dE}{dt}\prec E(t)&+\int_0^t E(s)\, ds+\int_0^t \Vert f(s)\Vert_{H^1}^2\, ds+\Vert f\Vert_{H^2}^2+\Vert f\Vert^2_{L^\infty\times L^2}\\
&+\Vert g_0\Vert^2_{H^3}+\Vert g_1\Vert^2_{H^3}+\Vert g_2\Vert^2_{H^2}
\end{split}
\]
and therefore 
\[
\begin{split}
\Vert W_3\Vert ^2_{L^2}\le E(t)&\prec E(0)+\int_0^t \Vert f(s)\Vert_{H^2}^2\, ds+ \Vert f\Vert_{L^\infty\times H^1}^2\\
&\prec \int_0^t \Vert f(s)\Vert_{H^2}^2\, ds+ \Vert f\Vert_{L^\infty\times H^1}^2+\Vert g_0\Vert^2_{H^4}+\Vert g_1\Vert^2_{H^3}+\Vert g_2\Vert^2_{H^2}.
\end{split}
\]
This leads to
\beq
\label{est_u_H_2}
\Vert u(t)\Vert_{H^2}^2\le c\biggl( \int_0^t \Vert f(s)\Vert_{H^2}^2\, ds+ \Vert f\Vert_{L^\infty\times H^1}^2+\Vert g_0\Vert^2_{H^4}+\Vert g_1\Vert^2_{H^3}+\Vert g_2\Vert^2_{H^2}\biggr).
\eeq
An iteration of this argument leads to Sobolev estimates of every order as in (iii). More precisely, by iterated derivation with respect to $x$ we get systems with $A$ as principal matrix, a combination of $A'$ and $B$ as matrix of order $0$ and a right-hand side given by a combination of derivatives of $A$ and $B$ and the solutions of the previous steps. For instance at step $k$ we get up to $k-1$ derivatives of $A$ and $B$, the solutions of the previous $k-1$ systems and $F^{(k-1)}$. To estimate this right-hand side we use the Levi conditions (on the term where $B'$ and the solution of the previous step appear)  and also the relations between the systems solutions, with application of the Minkowski's inequality as many times as needed.
\end{proof}
\begin{remark}
The following table summarising the strategy adopted in the proof above.
\begin{center}
\begin{tabular}{|c|c|c|c|}
\hline
System in & Hypotheses & $L^2$ estimates for & $u$ belongs to \\
\hline
$U$ & $a\ge 0$ and $(LC)$ & $U_3$ & $L^2$\\
\hline
$V$ & $a\ge 0$, $(LC)_1$ and $B^\infty$ coefficients & $V_3$ and $U_2$ & $H^1$\\
\hline
$W$ & $a\ge 0$, $(LC)_1$ and $B^\infty$ coefficients & $W_3$, $V_2$ and $U_1$ & $H^2$\\
\hline
$\cdots$ & $\cdots$ & $\cdots$ & $\cdots$\\
\hline
\end{tabular}
\end{center}
\end{remark}
\vspace{0.2cm}

Further lower order terms could be added to our equation leading to a bigger system size, however, this is not the main scope of our paper. Our scope is to extend our symmetrisation method to any space dimension.

\section{Third order hyperbolic equations in higher space dimension}
We now pass to investigate the previous third order equation in higher space dimension. This is a rather delicate topic because as we will see it is not immediate to construct a symmetriser and most importantly to have a nearly diagonal symmetriser.
 
\subsection{Our problem}
Let us consider equations of the type
\[
\partial_t^3u-\sum_{i=1}^n a_i(x)\partial_t\partial_{x_i}^2u+\sum_{i=1}^n b_{1,i}(x)\partial_{x_i}^2 u+\sum_{i=1}^n  b_{2,i}(x)\partial_t\partial_{x_i}u+b_{3,n}(x)\partial_t^2u=f(t,x),
\]
where $a_i\ge 0$ for all $i=1,\cdots,n$. Setting
\[
U=(\partial_{x_1}^2,\partial_{x_2}^2,\cdots,\partial_{x_n}^2,\partial_{x_1}\partial_t,\partial_{x_2}\partial_t,\cdots,\partial_{x_n}\partial_t,\partial_t^2)^Tu
\]
we can rewrite the equation above as 
\[
 \partial_t U=\sum_{k=1}^nA_k(x)\partial_{x_k}U+B(x)U+F, 
\]
where the matrices $A_k$ and $B$ have size $2n+1$ and $F$ is the $(2n+1)$-column with the (2n+1)th-entry equal to $f$ and all the others $0$.
 
The matrices $A_k$ have entries $(a_{k,ij})_{ij}$ as follows:
\[
\begin{split}
a_{k,ij}&=1,\quad\text{for $i=k$ and $j=n+k$},\\
a_{k,ij}&=1,\quad\text{for $i=n+k$ and $j=2n+1$},\\
a_{k,ij}&=a_k,\quad\text{for $i=2n+1$ and $j=n+k$},\\
a_{k,ij}&=0,\quad\text{otherwise}.
\end{split}
\]
This means that every matrix $A_k$ can be seen as a matrix in Sylvester form with $2n-2$ identically zero columns and rows.
In particular when $n=2$ we have
\[
A_1=\left(
	\begin{array}{ccccc}
	0& 0 & 1 & 0 & 0\\
	0& 0 & 0 & 0 & 0\\ 
	0& 0 & 0 & 0 & 1\\
	0& 0 & 0 & 0 & 0\\
	0 & 0 & a_1 & 0 & 0
	\end{array}
	\right),\quad A_2=\left(
	\begin{array}{ccccc}
	0& 0 & 0 & 0 & 0\\
	0& 0 & 0 & 1 & 0\\ 
	0& 0 & 0 & 0 & 0\\
	0& 0 & 0 & 0 & 1\\
	0 & 0 & 0 & a_2 & 0
	\end{array}
	\right).
\]
The matrix $B$ is has only the last row not identically zero, i.e.,
\[
(-b_{1,1},\cdots, -b_{1,n}, -b_{2,1},\cdots, -b_{2,n}, -b_{3,n}).
\]

\subsection{The symmetriser $Q$}

Our aim is to adapt the symmetriser method employed in space dimension 1 to this particular equation and to show that it will lead to Levi conditions analogous to the ones encountered when $n=1$. Inspired by \cite{G20} we define a common symmetriser $Q$ for all the matrices $A_k$, $k=1,\dots,n$.

\begin{proposition}
\label{prop_symm_3_n}
Let  
\[
Q=\frac{1}{3}\left(
	\begin{array}{cccccccc}
	a_1^2 & a_1a_2& \cdots &a_1a_{n} & 0& \cdots & \cdots & -a_1\\
	a_1a_2 & a_2^2 & a_2a_3 & \cdots & \cdots &  0 & \cdots & -a_2\\
	\vdots & \vdots & \vdots & \vdots & \vdots & \vdots & \vdots & \vdots\\
	a_1a_k & \vdots & a_k^2  & \vdots& 0 & \vdots  &\vdots& -a_k\\
	\vdots & \vdots & \vdots & \vdots & \vdots & \vdots & \vdots & \vdots\\
	a_1a_n & \vdots & \vdots & a_n^2 &  \vdots &\vdots & \vdots & -a_n\\
	\vdots & \vdots & \vdots & \vdots & 2a_1 & \vdots & \vdots & 0\\
	\vdots & \vdots & \vdots & \vdots & \vdots & \vdots & \vdots & \vdots\\
	\vdots & \vdots & \vdots & \vdots & \vdots & \vdots & 2a_n &  0\\
	-a_1 & -a_2 & \cdots & -a_n & 0 & \vdots & 0 & 3
	 \end{array}
	\right)
\]
\begin{itemize}
\item[(i)] $Q$ is a symmetriser of $A_k$ for every $k=1,\dots,n$, i.e, $QA_k=A_k^\ast Q$ is the symmetric matrix with entries $n+k, 2n+1$ and $2n+1, n+k$ equals to $\frac{2}{3}a_k$ and otherwise $0$.
\item[(ii)] For every $v\in \R^{2n+1}$
\[
\begin{split}
3\lara{Qv,v}&=\sum_{k=1}^n\lara{a_k^2v_k,v_k}+2\sum_{k=1}^n \lara{a_k v_{n+k},v_{n+k}}+3\Vert v_{2n+1}\Vert^2-2\sum_{k+1}^n\lara{a_kv_k, v_{2n+1}}\\
&+2\sum_{1\le i<j\le n}\lara{a_iv_i,a_{j}v_j}\\
&=\Vert \sum_{k=1}^n a_kv_k-v_{2n+1}\Vert^2+2\sum_{k=1}^n \lara{a_k v_{n+k},v_{n+k}}+2\Vert v_{2n+1}\Vert^2.
\end{split}
\]

\end{itemize}
\end{proposition}
\begin{proof}
Note that the matrix $Q$ has the following block structure 
\[
\frac{1}{3}\left(
	\begin{array}{ccc}
	S_{nn} & 0_{nn} & C\\
	0_{nn} & D & 0_{n1}\\
	C^T & 0_{1n} & 3 
	 \end{array}
	\right),
\]
where $S$ is the symmetric matrix 
\[
\left(
	\begin{array}{cccc}
	a_1^2 & a_1a_2 & \cdots & a_1a_n\\
	a_1a_2 & a_2^2 & \cdots &a_2a_n\\
	\vdots & \vdots & \vdots & \vdots\\
	a_1a_n & a_2a_n& \cdots & a_n^2
	 	 \end{array}
	\right),
\] 
$C=(-a_1,-a_2,\cdots,-a_n)^T$ and $D$ is the diagonal matrix with diagonal $(2a_1, 2a_2,\cdots,a_n)$.
\begin{itemize}
\item[(i)]
By definition of the matrices $Q$ and $A_k$ we have 
\[
\frac{1}{3}\left(
	\begin{array}{ccc}
	S & 0_{nn} & C\\
	0_{nn} & D & 0_{n1}\\
	C^T & 0_{12} & 3 
	 \end{array}
	\right)A_k=\frac{1}{3}\left(
	\begin{array}{ccc}
	S & 0_{nn} & C\\
	0_{nn} & D & 0_{n1}\\
	C^T & 0_{1n} & 3 
	 \end{array}
	\right)\left(
	\begin{array}{ccc}
	0_{nn} & 1_{k,n+k} & 0\\
	0_{nn} & 0_{nn} & 1_{n+k,2n+1}\\
	0& a_{k,n+k} & 3 
	 \end{array}
	\right),
\]
where $1_{k,n+k}$ denotes the $n\times n$-matrix with entry $k,n+k$ equals to 1 and $0$ otherwise, $1_{n+k,2n+1}$ is the column with $1$ in position $n+k,2n+1$ and otherwise $0$ and finally $a_{k,n+k}$ is the row with $a_k$ in position $n+k$ and otherwise $0$. By direct computations we see that when multiplying the row $(S\,  0_{nn}\,  C)$ by the column $(1_{k,n+k}\, 0_{nn}\, a_{k,n+k})$ we get all zero entries since $a_k^2-a_k^2=0$ and $a_ia_k-a_ka_i=0$ for $i\neq k$. Similarly, by multiplying the other rows for the columns of $A_k$ we get
\[
QA_k=\frac{1}{3}\left(
	\begin{array}{ccc}
	0_{nn} & 0_{nn} & 0\\
	0_{nn} & 0_{nn} & C_k\\
	0& C_k^T & 0
	 \end{array}
	\right),
\]
where $C_k$ has all entries zero a part from the one in position $n+k, 2n+1$ which is equal to $2a_k$.
\item[(ii)] We easily see that 
\[
\begin{split}
&3\lara{Qv,v}=\langle\left(
	\begin{array}{ccc}
	S & 0_{nn} & C\\
	0_{nn} & D & 0_{n1}\\
	C^T & 0_{1n} & 3 
	 \end{array}
	\right)v,v\rangle\\
	&=\sum_{k=1}^n\lara{a_k^2v_k,v_k}+2\sum_{k=1}^n \lara{a_k v_{n+k},v_{n+k}}+3\Vert v_{2n+1}\Vert^2-2\sum_{k=1}^n\lara{a_kv_k, v_{2n+1}}\\
&+2\sum_{1\le i<j\le n}\lara{a_iv_i,a_{j}v_j}\\
&=\Vert \sum_{k=1}^n a_kv_k-v_{2n+1}\Vert^2+2\sum_{k=1}^n \lara{a_k v_{n+k},v_{n+k}}+2\Vert v_{2n+1}\Vert^2.
\end{split}
\]
\end{itemize}
\end{proof}

\begin{remark}
Note that the symmetriser $Q$ is positive semi-definite since $\lara{Qv,v}\ge 0$ but ${\rm det}Q=0$. In addition we cannot apply directly the results on the standard symmetriser because the matrices $A_k$ are not in Sylvester form in the classical sense. What is clear from the definition of $Q$ above is that if we define the energy $E(t)$ as $(QU, U)_{L^2}$ then $E(t)\ge 0$ and 
\[
E(t)\ge \frac{2}{3}\Vert U_{2n+1}\Vert_{L^2}^2=\frac{2}{3}\Vert \partial_t u\Vert_{L^2}^2.
\]
In addiiton, assuming that the equation coefficients are bounded we have that 
\[
E(t)\prec \sum_{i=1}^{2n+1}\Vert U_i\Vert_{L^2}^2.
\]
\end{remark}
\subsection{The energy estimates}
Let us now focus on the Cauchy problem 
\[
\begin{split}
\partial_t^3u-\sum_{i=1}^n a_i(x)\partial_t\partial_{x_i}^2u+\sum_{i=1}^n b_{1,i}(x)\partial_{x_i}^2 u+\sum_{i=1}^n  b_{2,i}(x)\partial_t\partial_{x_i}u+b_{3,n}(x)\partial_t^2u&=f(t,x),\\
u(0,x)&=g_0(x),\\
\partial_tu(0,x)&=g_1(x),\\
\partial_t^2u(0,x)&=g_2(x),
\end{split}
\]
which is transformed into 
\beq
\label{CP_3_n}
\begin{split}
\partial_t U&=\sum_{k=1}^nA_k(x)\partial_{x_k}U+B(x)U+F,\\
U(0)&=(\partial^2_{x_1}g_0, \cdots, \partial^2_{x_n}g_0, \partial_{x_1}g_1,\cdots, \partial_{x_n}g_1, g_2).
\end{split}
\eeq
Given the energy $E(t)=(QU, U)_{L^2}$, we have
 \[
 \begin{split}
 &\frac{dE(t)}{dt}=(\partial_t(QU),U)_{L^2}+(QU,\partial_tU)_{L^2}\\
 &=-\sum_{k=1}^n(\partial_{x_k}(QA_k)U,U)_{L^2}+((QB+B^\ast Q)U,U)_{L^2}+2(QU,F)_{L^2}.
\end{split}
 \]
 \begin{trivlist}
 \item[(i)] {\bf Estimate of $(\partial_{x_k}(QA_k)U,U)_{L^2}$.}
 By direct computations
 \[
 (\partial_{x_k}(QA_k)U,U)_{L^2}=\frac{4}{3}(\partial_{x_k}a_kU_{n+k}, U_{2n+1})_{L^2}.
 \]
 By Glaeser's inequality we have 
 \[
  (\partial_{x_k}(QA_k)U,U)_{L^2}=\frac{4}{3}(\partial_{x_k}a_kU_{n+k}, U_{2n+1})_{L^2}\prec (a_k U_{n+k}, U_{n+k})_{L^2}+\Vert U_{2n+1}\Vert_{L^2}^2\prec E(t).
 \]
 So,
 \[
 -\sum_{k=1}^n(\partial_{x_k}(QA_k)U,U)_{L^2}\prec E(t).
 \]
\item[(i)] {\bf Estimate of $((QB+B^\ast Q)U,U)_{L^2}$.} It is difficult to estimate this term without the nearly-diagonality of $Q$. However, we can overcome this issue by imposing ad-hoc Levi conditions on the lower order terms. For the sake of simplicity and to explain better
our method we work in space dimension 2. By definition of the symmetriser we have that 
\[
\begin{split}
3E(t)&=\Vert a_1U_1+a_2U_2-U_5\Vert_{L^2}^2+(2a_1U_3,U_3)_{L^2}+(2a_2U_4,U_4)_{L^2}+2\Vert U_5\Vert_{L^2}^2\\
&\ge \frac{1}{2} \Vert a_1U_1+a_2U_2\Vert_{L^2}^2-\Vert U_5\Vert_{L^2}^2+(2a_1U_3,U_3)_{L^2}+(2a_2U_4,U_4)_{L^2}+2\Vert U_5\Vert_{L^2}^2\\
&= \frac{1}{2} \Vert a_1U_1+a_2U_2\Vert_{L^2}^2+(2a_1U_3,U_3)_{L^2}+(2a_2U_4,U_4)_{L^2}+\Vert U_5\Vert_{L^2}^2.
\end{split}
\]
This show that we can estimate the energy from below with the quantity
\[
 \Vert a_1U_1+a_2U_2\Vert_{L^2}^2+(a_1U_3,U_3)_{L^2}+(a_2U_4,U_4)_{L^2}+\Vert U_5\Vert_{L^2}^2.
\]
Let now 
\[
B=\left(
	\begin{array}{cccccc}
	0 & 0 & 0 & 0 & 0\\
	0&  0 & 0 &  0 & 0\\
	0& 0 & 0 & 0 & 0  \\
	0& 0 & 0 & 0 & 0  \\
	-b_1& -b_2 & -b_3 & -b_4 & -b_5
	 \end{array}
	\right).
\]
By straightforward computations we have 
\[
3QB=\left(
	\begin{array}{cccccc}
	a_1b_1 & a_1b_2& a_1b_3 & a_1b_4 & a_1b_5\\
	a_2b_1 & a_2b_2& a_2b_3 & a_2b_4 & a_2b_5\\\
	0& 0 & 0 & 0 & 0  \\
	0& 0 & 0 & 0 & 0  \\
	-3b_1& -3b_2 & -3b_3 & -3b_4 & -3b_5
	 \end{array}
	\right)
\]
and therefore
\[
3(QB+B^\ast Q)=\left(
	\begin{array}{cccccc}
	2a_1b_1 & a_1b_2+a_2b_1& a_1b_3 & a_1b_4 & a_1b_5-3b_1\\
	a_1b_2+a_2b_1 & 2a_2b_2& a_2b_3 & a_2b_4 & a_2b_5-3b_2\\\
	a_1b_3& a_2b_3 & 0 & 0 & -3b_3  \\
	a_1b_4& a_2b_4 & 0 & 0 & -3b_4  \\
	a_1b_5-3b_1& a_2b_5-3b_2 & -3b_3 & -3b_4 & -6b_5
	 \end{array}
	\right).
\]
It follows that 
\[
\begin{split}
3((QB+B^\ast Q)U,U)_{L^2}&=(2a_1b_1U_1,U_1)_{L^2}+2((a_1b_2+a_2b_1)U_1,U_2)_{L^2}+(2a_2b_2U_2,U_2)_{L^2}\\
&+2(a_1b_3U_3,U_1)_{L^2}+2(a_2b_3U_3,U_2)_{L^2}\\
&+2(a_1b_4U_4,U_1)_{L^2}+2(a_2b_4U_4,U_2)_{L^2}\\
&+2((a_1b_5-3b_1)U_1,U_5)_{L^2}+2((a_2b_5-3b_2)U_2,U_5)_{L^2}\\
&+2(-3b_3U_3,U_5)_{L^2}+2(-3b_4U_4,U_5)_{L^2}\\
&+(-6b_5U_5,U_5)_{L^2}.
\end{split}
\]
Analysing the term 
\[
I_1=(2a_1b_1U_1,U_1)_{L^2}+2((a_1b_2+a_2b_1)U_1,U_2)_{L^2}+(2a_2b_2U_2,U_2)_{L^2}
\]
we see that if 
\[
\begin{split}
b_1&=\lambda a_1,\\
b_2&=\lambda a_2,\\
\end{split}
\]
for some bounded function $\lambda$ then 
\[
I_1\prec \Vert |\lambda| a_1U_1+|\lambda| a_2U_2\Vert_{L^2}^2\prec \Vert a_1U_1+a_2U_2\Vert_{L^2}^2.
\]
We now write 
\[
I_2=2(a_1b_3U_3,U_1)_{L^2}+2(a_2b_3U_3,U_2)_{L^2}
\]
as 
\[
I_2=2(b_3U_3,a_1U_1+a_2U_2).
\]
Hence, if
\[
|b_3|\prec \sqrt{a_1}
\]
we have that 
\[
I_2\prec  \Vert a_1U_1+a_2U_2\Vert_{L^2}^2+(a_1U_3,U_3)_{L^2}.
\]
Analogously, if 
\[
|b_4|\prec \sqrt{a_2}
\]
then
\[
\begin{split}
I_3=2(a_1b_4U_4,U_1)_{L^2}+2(a_2b_4U_4,U_2)_{L^2}&=2(b_4U_4, a_1U_1+a_2U_2)_{L^2}\\
&\prec  \Vert a_1U_1+a_2U_2\Vert_{L^2}^2+(a_2U_4,U_4)_{L^2}.
\end{split}
\]
We now write 
\[
I_4=2((a_1b_5-3b_1)U_1,U_5)_{L^2}+2((a_2b_5-3b_2)U_2,U_5)_{L^2}
\]
as 
\[
I_4=2(b_5U_5, a_1U_1+a_2U_2)_{L^2}+2(-3b_1U_1-3b_2U_2,U_5)_{L^2}.
\]
By the Levi conditions $b_1=\lambda a_1$ and $b_2=\lambda a_2$ we have that 
\[
2(-3b_1U_1-3b_2U_2,U_5)_{L^2}=-6 (\lambda a_1U_1+\lambda a_2U_2,U_5)_{L^2}.
\]
Thus, if 
\[
|b_5|\prec 1
\]
we conclude that 
\[
I_4\prec \Vert a_1U_1+a_2U_2\Vert_{L^2}^2+\Vert U_5\Vert^2_{L^2}.
\]
Finally, from the Levi conditions $|b_3|\prec \sqrt{a_1}$, $|b_4|\prec \sqrt{a_2}$ and $|b_5|\prec 1$ we easily obtain that 
\[
2(-3b_3U_3,U_5)_{L^2}+2(-3b_4U_4,U_5)_{L^2}+(-6b_5U_5,U_5)_{L^2}\prec (a_1U_3,U_3)_{L^2}+(a_2U_4,U_4)_{L^2}+\Vert U_5\Vert_{L^2}^2.
\]
Summarising, we have proven that under the Levi conditions (LC),
\[
\begin{split}
b_1&=\lambda a_1,\\
b_2&=\lambda a_2,\\
|b_3|&\prec \sqrt{a_1},\\
|b_4|&\prec\sqrt{a_2},\\
|b_5|&\prec 1,
\end{split}
\]
the estimate 
\[
\begin{split}
3((QB+b^\ast Q)U,U)_{L^2}&\prec  \Vert a_1U_1+a_2U_2\Vert_{L^2}^2+(a_1U_3,U_3)_{L^2}+(a_2U_4,U_4)_{L^2}+\Vert U_5\Vert_{L^2}^2\\
&\prec E(t)
\end{split}
\]
holds, i.e., the matrix of the lower order terms can be bounded by the energy $E(t)$.
\item[(iii)] {\bf Estimate of $(QU,F)_{L^2}$}
 By direct computations we have that 
 \[
 3(QU,F)_{L^2}=(-a_1U_1-a_2U_2+3U_5,f)_{L^2}\prec \Vert a_1U_1+a_2U_2\Vert_{L^2}^2+\Vert U_5\Vert_{L^2}^2+\Vert f\Vert_{L^2}^2\prec E(t)+\Vert f\Vert_{L^2}^2.
 \]
 \item[(iv)]{\bf Conclusion:} Under the Levi conditions (LC) we have that 
 \[
 \frac{dE(t)}{dt}\prec E(t)+\Vert f\Vert_{L^2}^2.
 \]
By application of Gr\"onwall's lemma and the boundedness of the equation coefficients we get, for some constant $c>0$, the estimate
\[
\begin{split}
\Vert \partial_t^2u(t)\Vert_{L^2}^2&\le E(t)\le c\biggl( E(0)+\int_0^t \Vert f(s)\Vert_{L^2}^2\, ds\biggr)\\
&\prec\Vert g_0\Vert^2_{H^2}+\Vert g_1\Vert_{H^1}^2+\Vert g_2\Vert_{L^2}^2+\int_0^t \Vert f(s)\Vert_{L^2}^2\, ds.
\end{split}
\]
By the fundamental theorem of calculus and Minkowski's integral inequality we get  
\[
 \Vert u(t)\Vert_{L^2}^2\le 2\Vert u(t)-u(0)\Vert^2_{L^2}+2\Vert u(0)\Vert_{L^2}^2\le 2\int_0^t \Vert \partial_t u(s)\Vert_{L^2}^2\, ds+ 2\Vert u(0)\Vert_{L^2}^2
\]
and 
\[
 \Vert \partial_tu(s)\Vert_{L^2}^2\le 2\Vert \partial_tu(s)-\partial_tu(0)\Vert^2_{L^2}+2\Vert \partial_tu(0)\Vert_{L^2}^2\le 2\int_0^s \Vert \partial^2_t u(r)\Vert_{L^2}^2\, dr+ 2\Vert \partial_tu(0)\Vert_{L^2}^2.
\]

It follows that 
\[
\begin{split}
 \Vert u(t)\Vert_{L^2}^2&\prec \int_0^t\int_0^s\Vert \partial_t^2u(r)\Vert_{L^2}^2\, dr\, ds+\Vert u(0)\Vert_{L^2}^2+\Vert \partial_t u(0)\Vert_{L^2}^2\\
 &\prec \int_0^t\int_0^s\Vert \partial_t^2u(r)\Vert_{L^2}^2\, dr\, ds+\Vert g_0\Vert_{L^2}^2+\Vert g_1\Vert_{L^2}^2\\
 &\prec \Vert g_0\Vert^2_{H^2}+\Vert g_1\Vert_{H^1}^2+\Vert g_2\Vert_{L^2}^2+\int_0^t \Vert f(s)\Vert_{L^2}^2\, ds.
 \end{split}
\]
\end{trivlist}
This shows that our Cauchy problem is $L^2$ well-posed. 
\subsection{$L^2$ well-posedness} It is only a technical matter to extend the argument above to any space dimension by employing matrices with a bigger size but the same proof strategy. We therefore have the following theorem. 

%
%
%
 
\begin{theorem}
\label{theo_2_3_L_n}
Let 
\[
\begin{split}
\partial_t^3u-\sum_{i=1}^n a_i(x)\partial_t\partial_{x_i}^2u+\sum_{i=1}^n b_i(x)\partial_{x_i}^2 u+\sum_{i=1}^n  b_{2,i}(x)\partial_t\partial_{x_i}u+b_{3,n}(x)\partial_t^2u&=f(t,x),\\
u(0,x)&=g_0(x),\\
\partial_tu(0,x)&=g_1(x),\\
\partial_t^2u(0,x)&=g_2(x),
\end{split}
\]
where the equation coefficients are real-valued, smooth and with bounded derivatives of any order, $a_i\ge 0$ for $i=1,\cdots, n$ and $f\in C([0,T],L^2(\R^2))$. Under the Levi conditions (LC),
\[
\begin{split}
b_i&=\lambda a_i,\qquad i=1,\cdots, n,\\
|b_{2,i}|&\prec \sqrt{a_i},\qquad i=1,\cdots, n,\\
|b_{3,n}|&\prec 1,
\end{split}
\]
where $\lambda$ is a bounded function, the Cauchy problem has a unique solution in $C^3([0,T], L^2(\R^n))$ provided that $g_0\in H^2(\R^n)$, $g_1\in H^1(\R^n)$ and $g_2\in L^2(\R^n)$.
\end{theorem}

\subsection{$H^1$ well-posedness}
To find out under which assumptions and Levi conditions the Cauchy problem 
\[
\begin{split}
\partial_t^3u-\sum_{i=1}^n a_i(x)\partial_t\partial_{x_i}^2u+\sum_{i=1}^n b_i(x)\partial_{x_i}^2 u+\sum_{i=1}^n  b_{2,i}(x)\partial_t\partial_{x_i}u+b_{3,n}(x)\partial_t^2u&=f(t,x),\\
u(0,x)&=g_0(x),\\
\partial_tu(0,x)&=g_1(x),\\
\partial_t^2u(0,x)&=g_2(x),
\end{split}
\]
is $H^1$ well-posed we will need to work not on the system 
\[
\begin{split}
\partial_t U&=\sum_{k=1}^nA_k(x)\partial_{x_i}U+B(x)U+F,\\
U(0)&=(\partial^2_{x_1}g_0, \cdots, \partial^2_{x_n}g_0, \partial_{x_1}g_1,\cdots, \partial_{x_n}g_1, g_2)^T
\end{split}
\]
but rather on the $n(2n+1)\times n(2n+1)$-system in $V$ obtained by deriving $U$ with respect to $x_i$, with $i=1,\cdots,n$. In detail, 
\[
\partial_t V=\sum_{i=1}^n \wt{A_i}(x)\partial_{x_i}V+ \wt{B}V+\wt{F},
\]
where $\wt{A_i}$ is a diagonal matrix with $n$ repeated blocks $A_i$ on the diagonal, 
\[
{\wt{B}}=\left(
	\begin{array}{ccccc}
	\partial_{x_1}{A_1}+{B} & \partial_{x_1}{A_2} & \cdots & \cdots & \partial_{x_1}{A_n}\\
	\partial_{x_2}{A_1}& \partial_{x_2}{A_2}+{B} & \cdots & \cdots & \partial_{x_2}{A_n}\\
	\vdots & \vdots & \vdots & \vdots & \vdots\\
	\partial_{x_k}{A_1}&  \cdots &\partial_{x_k}{A_k}+{B} & \cdots & \partial_{x_k}{A_n}\\
	\vdots & \vdots & \vdots & \vdots & \vdots\\
	\partial_{x_n}{A_1}& \partial_{x_n}{A_2} & \cdots & \cdots & \partial_{x_n}{A_n}+{B}\\
	 
	\end{array}
	\right)
\]
and 
\[
{\wt{F}}=\nabla_x {F}+\left(
	\begin{array}{c}
	(\partial_{x_1}{B})U\\ 
	(\partial_{x_2}{B})U\\ 
	\vdots \\
	(\partial_{x_k}{B})U\\
	\vdots\\
	(\partial_{x_n}{B})U\\ 
	\end{array}
	\right)= \left(
	\begin{array}{c}
	\partial_{x_1}F\\ 
	\partial_{x_2}F\\ 
	\vdots \\
	\partial_{x_k}F\\
	\vdots\\
	\partial_{x_n}F\\ 
	\end{array}
	\right)+\left(
	\begin{array}{c}
	(\partial_{x_1}{B})U\\ 
	(\partial_{x_2}{B})U\\ 
	\vdots \\
	(\partial_{x_k}{B})U\\
	\vdots\\
	(\partial_{x_n}{B})U\\ 
	\end{array}
	\right).
\]
Arguing as in \cite{G20} Subsection 4.3, we make use of the energy $E(t)=(\wt{Q}V,V)_{L^2}$, where ${\wt{Q}}$ is a block-diagonal matrix with $n$ identical blocks equal to $Q$. Arguing as for the system in $U$ we get 
\[
 \frac{dE(t)}{dt} =-\sum_{k=1}^n(\partial_{x_k}({\wt{Q}}{\wt{A_k}})V,V)_{L^2}+(({\wt{Q}}{\wt{B}}+{\wt{B}}^\ast{\wt{Q}})V,V)_{L^2}+2({\wt{Q}}V,{\wt{F}})_{L^2}.
\]
Because of the block-diagonal structure of both $\wt{A_k}$ and $\wt{Q}$ it follows immediately that under our hypotheses on the coefficients $a_i$ we have that 
\[
-\sum_{k=1}^n(\partial_{x_k}({\wt{Q}}{\wt{A_k}})V,V)_{L^2}\prec E(t).
\]
In order to estimate $(({\wt{Q}}{\wt{B}}+{\wt{B}}^\ast{\wt{Q}})V,V)_{L^2}$ it is sufficient to investigate the structure of $({\wt{Q}}{\wt{B}}V,V)_{L^2}$. This can be written as $(S_1V,V)_{L^2}+(S_2V,V)_{L^2}$, where  
\[
S_1=\wt{Q}\left(
	\begin{array}{ccccc}
	\partial_{x_1}{A_1} & \partial_{x_1}{A_2} & \cdots & \cdots & \partial_{x_1}{A_n}\\
	\partial_{x_2}{A_1}& \partial_{x_2}{A_2}& \cdots & \cdots & \partial_{x_2}{A_n}\\
	\vdots & \vdots & \vdots & \vdots & \vdots\\
	\partial_{x_k}{A_1}&  \cdots &\partial_{x_k}{A_k}& \cdots & \partial_{x_k}{A_n}\\
	\vdots & \vdots & \vdots & \vdots & \vdots\\
	\partial_{x_n}{A_1}& \partial_{x_n}{A_2} & \cdots & \cdots & \partial_{x_n}{A_n}\\
	 \end{array}
	\right)
\]
and $S_2$ is a block diagonal matrix with repeated blocks $QB$. It follows that under the Levi conditions (LC) we get immediately that $((S_2+S_2^\ast)V,V)_{L^2}\prec E(t)$. It remains to estimate 
\[
((S_1+S_1^\ast) V, V)_{L^2}.
\]
Let us argue for simplicity in the case $n=2$ (dimensions higher than 2 are technically more challenging but do not change the nature of the argument). We also denote the lower order terms coefficients as $b_i$, with $i=1,\cdots,5$ since we deal with $5\times 5$ matrices. We have 
\[
\begin{split}
((S_1+S_1^\ast) V, V)_{L^2}
 &=2\biggr(\left(
	\begin{array}{cc}
	Q& 0\\  
	0 & Q\\
	 \end{array}
	\right)
	\left(
	\begin{array}{cc}
	\partial_{x_1}{A_1} & \partial_{x_1}{A_2}\\  
	\partial_{x_2}{A_1}& \partial_{x_2}{A_2}\\
	 \end{array}
	\right)V,V\biggl)_{L^2}\\
	& =2\biggr( 
	\left(
	\begin{array}{cc}
	Q\partial_{x_1}{A_1} & Q\partial_{x_1}{A_2}\\  
	Q\partial_{x_2}{A_1}& Q\partial_{x_2}{A_2}\\
	 \end{array}
	\right)V,V\biggl)_{L^2}.
 \end{split}
\]
Note that by definition of the matrices $\partial_{x_i}A_j$ and $Q$ with $i,j=1,2$ we have 
\[
3Q\partial_{x_i}A_1=\left(
	\begin{array}{ccccc}
	0 & 0 & -a_1\partial_{x_i}a_1 & 0 & 0\\
	0 & 0 & -a_2\partial_{x_i}a_1 & 0 & 0\\
	0 & 0 & 0 & 0 & 0\\
	0 & 0 & 0 & 0 & 0\\
	0 & 0 & 3\partial_{x_i}a_1 & 0
	\end{array}
	\right), \qquad 
3Q\partial_{x_i}A_2=\left(
	\begin{array}{ccccc}
	0 & 0 &  0 & -a_1\partial_{x_i}a_2 & 0 \\
	0 & 0 & 0 &  -a_2\partial_{x_i}a_2 & 0 \\
	0 & 0 & 0 & 0 & 0\\
	0 & 0 & 0 & 0 & 0\\
	0 & 0 &  0 & 3\partial_{x_i}a_2 & 0
	\end{array}
	\right).
	\]
Hence, the only non-identically zero columns in the matrix  
\[
\left(
	\begin{array}{cc}
	Q\partial_{x_1}{A_1} & Q\partial_{x_1}{A_2}\\  
	Q\partial_{x_2}{A_1}& Q\partial_{x_2}{A_2}\\
	 \end{array}
	\right)
	\]
are the third and the ninth column. It follows that the matrix above multiplied by $V$ gives the column vector 
\[
\left(
	\begin{array}{c}
- \frac{1}{3}a_1\partial_{x_1}a_1V_3-\frac{1}{3}a_1\partial_{x_1}a_2V_9\\
- \frac{1}{3}a_2\partial_{x_1}a_1V_3-\frac{1}{3}a_2\partial_{x_1}a_2V_9\\
0\\
0\\
\partial_{x_1}a_1V_3+\partial_{x_1}a_2V_9\\
- \frac{1}{3}a_1\partial_{x_2}a_1V_3-\frac{1}{3}a_1\partial_{x_2}a_2V_9\\
- \frac{1}{3}a_2\partial_{x_2}a_1V_3-\frac{1}{3}a_2\partial_{x_2}a_2V_9\\
0\\
0\\
\partial_{x_2}a_1V_3+\partial_{x_2}a_2V_9
	 \end{array}
	\right).
\]
Concluding,
\[
\begin{split}
\frac{1}{2}((S_1+S_1^\ast) V, V)_{L^2}
 &=(-\frac{1}{3}a_1\partial_{x_1}a_1V_3-\frac{1}{3}a_1\partial_{x_1}a_2V_9, V_1)_{L^2}\\
 &+(-\frac{1}{3}a_2\partial_{x_1}a_1V_3-\frac{1}{3}a_2\partial_{x_1}a_2V_9, V_2)_{L^2}\\
 &+(\partial_{x_1}a_1V_3+\partial_{x_1}a_2V_9, V_5)_{L^2}\\
 &+(- \frac{1}{3}a_1\partial_{x_2}a_1V_3-\frac{1}{3}a_1\partial_{x_2}a_2V_9, V_6)_{L^2}\\
 &+(- \frac{1}{3}a_2\partial_{x_2}a_1V_3-\frac{1}{3}a_2\partial_{x_2}a_2V_9, V_7)_{L^2}\\
 &+(\partial_{x_2}a_1V_3+\partial_{x_2}a_2V_9, V_{10})_{L^2}.
 \end{split}
\]
By Glaeser's inequality, we have that $|\partial_{x_j}a_i|\prec \sqrt{a_i}$ for all $i,j=1,2$. It follows that 
\[
\begin{split}
(-\frac{1}{3}a_1\partial_{x_1}a_1V_3-\frac{1}{3}a_1\partial_{x_1}a_2V_9, V_1)_{L^2}&\prec (a_1^2V_1, V_1)_{L^2}+(2a_1V_3,V_3)_{L^2}+(2a_2V_9,V_9)_{L^2}\\
 (-\frac{1}{3}a_2\partial_{x_1}a_1V_3-\frac{1}{3}a_2\partial_{x_1}a_2V_9, V_2)_{L^2}&\prec (a_2^2V_2, V_2)_{L^2}+(2a_1V_3,V_3)_{L^2}+(2a_2V_9,V_9)_{L^2}\\
 (\partial_{x_1}a_1V_3+\partial_{x_1}a_2V_9, V_5)_{L^2}&\prec(2a_1V_3,V_3)_{L^2}+(2a_2V_9,V_9)_{L^2}+\Vert V_5\Vert_{L^2}^2\\
 (- \frac{1}{3}a_1\partial_{x_2}a_1V_3-\frac{1}{3}a_1\partial_{x_2}a_2V_9, V_6)_{L^2}&\prec (a_1^2V_6, V_6)_{L^2}+(2a_1V_3,V_3)_{L^2}+(2a_2V_9,V_9)_{L^2}\\
 (- \frac{1}{3}a_2\partial_{x_2}a_1V_3-\frac{1}{3}a_2\partial_{x_2}a_2V_9, V_7)_{L^2}&\prec (a_2^2V_7, V_7)_{L^2}+(2a_1V_3,V_3)_{L^2}+(2a_2V_9,V_9)_{L^2}\\
 (\partial_{x_2}a_1V_3+\partial_{x_2}a_2V_9, V_{10})_{L^2}&\prec (2a_1V_3,V_3)_{L^2}+(2a_2V_9,V_9)_{L^2}+\Vert V_{10}\Vert_{L^2}^2.
 \end{split}
\]
Since all the terms in the right-hand side above are estimated by the energy $E(t)$ we have that 
\[
((S_1+S_1^\ast) V, V)_{L^2}\prec E(t).
\]
It remains to estimate $(\wt{Q}V, \wt{F})_{L^2}$. We begin by writing $(\wt{Q}V, \wt{F})_{L^2}$ as 
\[
\biggl(\left(
	\begin{array}{cc}
	Q & 0\\  
	0& Q \\
	 \end{array}
	\right)V, \left(\begin{array}{c}\partial_{x_1}F\\
	\partial_{x_2}F
	\end{array}
	\right)\biggr)_{L^2}+\biggl(\left(
	\begin{array}{cc}
	Q & 0\\  
	0& Q \\
	 \end{array}
	\right)V, \left(\begin{array}{c}(\partial_{x_1}B)U\\
	(\partial_{x_2}B)U
	\end{array}
	\right)\biggr)_{L^2}.
\]
Since the first addendum is
\[
(V_5, \partial_{x_1}f)_{L^2}+(V_{10}, \partial_{x_1}f)_{L^2}
\]
we easily conclude that 
\[
\biggl(\left(
	\begin{array}{cc}
	Q & 0\\  
	0& Q \\
	 \end{array}
	\right)V, \left(\begin{array}{c}\partial_{x_1}F\\
	\partial_{x_2}F
	\end{array}
	\right)\biggr)_{L^2}\prec E(t)+\Vert f\Vert_{H^1}^2.
\]
We now focus on 
\[
\biggl(\left(
	\begin{array}{cc}
	Q & 0\\  
	0& Q \\
	 \end{array}
	\right)V, \left(\begin{array}{c}(\partial_{x_1}B)U\\
	(\partial_{x_2}B)U
	\end{array}
	\right)\biggr)_{L^2}.
\]
We have that it coincides with 
\[
\begin{split}
&(V_5, -\partial_{x_1}b_1U_1-\partial_{x_1}b_2U_2-\partial_{x_1}b_3U_3-\partial_{x_1}b_4U_4-\partial_{x_1}b_5U_5)_{L^2}\\
&+(V_{10}, -\partial_{x_2}b_1U_1-\partial_{x_2}b_2U_2-\partial_{x_2}b_3U_3-\partial_{x_2}b_4U_4-\partial_{x_2}b_5U_5)_{L^2}.
\end{split}
\]
Estimating this final term means to estimate the norms $\Vert \partial_{x_j}b_kU_k\Vert^2_{L^2}$, where $j=1,2$ and $k=1,\cdots,5$. We distinguish between three cases:
\begin{enumerate}
\item[(1)] $k=1,3$;
\item[(2)] $k=2,4$;
\item[(3)] $k=5$.
\end{enumerate}
{\bf Case (1)}
In the first case we employ the argument already used in \cite{G20} for the wave equation, i.e., 
\begin{multline}
 \label{est_preparatory}
\Vert \partial_{x_j}b_kU_k\Vert^2_{L^2}=\biggl\Vert\int_0^t \partial_{x_j}b_k\partial_tU_k(s)\, ds+\partial_{x_j}b_kU_k(0)\biggr\Vert^2\\
\le 2\biggl\Vert \int_0^t \partial_{x_j}b_k\partial_tU_k(s)\, ds\biggr\Vert^2_{L^2}+2\Vert \partial_{x_j}b_kU_k(0)\Vert_{L^2}^2\\
\le2\biggl(\int_0^t\Vert \partial_{x_j}b_kV_{k+n}(s)\Vert_{L^2}\, ds\biggr)^2+2\Vert\partial_{x_j} b_kU_k(0)\Vert_{L^2}^2\\
\le 2T\int_0^t((\partial_{x_j}b_k)^2V_{k+n},V_{k+n})_{L^2}\, ds+2\Vert  \partial_{x_j}b_kU_k(0)\Vert_{L^2}^2\\
\le 2T\int_0^t((\partial_{x_j}b_k)^2V_{k+n},V_{k+n})_{L^2}\, ds+2\Vert \partial_{x_j}b_k\Vert_{L^\infty}^2\Vert U_k(0)\Vert_{L^2}^2
\end{multline}
where in this case $n=2$ and $k=1,3$. Under the Levi conditions (LC)
\[
b_1=\lambda a_1,\, b_2=\lambda a_2,\, |b_3|\prec \sqrt{a_1},\, |b_4|\prec\sqrt{a_2},\, |b_5|\prec 1,
\]
and the assumptions that the coefficients belong to $B^\infty$, we conclude that 
\[
\begin{split}
((\partial_{x_j}b_1)^2V_{3},V_{3})_{L^2}&=((\lambda\partial_{x_j}a_1+(\partial_{x_j}\lambda) a_1)^2V_{3},V_{3})_{L^2}\prec (2a_1V_3,V_3)_{L^2}\\
((\partial_{x_j}b_3)^2V_{5},V_{5})_{L^2}&\prec \Vert V_5\Vert^2_{L^2}.
\end{split}
\]
It follows that 
\beq
\label{est_1_3}
\Vert \partial_{x_j}b_kU_k\Vert^2_{L^2}\prec \int_0^t E(s)\, ds+\Vert U_k(0)\Vert^2_{L^2},
\eeq
for $k=1,3$.\\
{\bf Case (2)}\\
In analogy to Case (1) we get
\[
\Vert \partial_{x_j}b_{k}U_k\Vert^2_{L^2}\prec \int_0^t((\partial_{x_j}b_k)^2V_{2(2n+1)-(n-\frac{k}{2})},V_{2(2n+1)-(n-\frac{k}{2})})_{L^2}\, ds+2\Vert \partial_{x_j}b_k\Vert_{L^\infty}^2\Vert U_k(0)\Vert_{L^2}^2,
\]
where $k=2,4$ and $n=2$. Under the Levi conditions (LC) and the assumptions that the coefficients belong to $B^\infty$, we have
\[
\begin{split}
((\partial_{x_j}b_2)^2V_9,V_9)_{L^2}&=((\lambda \partial_{x_j}a_2+(\partial_{x_j}\lambda)a_2)^2V_9,V_9)_{L^2}\prec (2a_2V_9,V_9)_{L^2}\\
((\partial_{x_j}b_4)^2V_{10},V_{10})_{L^2}&\prec \Vert V_{10}\Vert^2.
\end{split}
\]
This leads to the estimate \eqref{est_1_3} also for $k=2,4$.\\
{\bf Case (3)}\\
When $k=5$ since $U_5$ has been already estimated in the $L^2$ well-posedness proof we immediately have that 
\beq
\label{est_k_5}
\Vert \partial_{x_j}b_5U_5\Vert^2_{L^2}\prec\Vert g_0\Vert^2_{H^2}+\Vert g_1\Vert_{H^1}^2+\Vert g_2\Vert_{L^2}^2+\int_0^t \Vert f(s)\Vert_{L^2}^2\, ds.
\eeq
Combining \eqref{est_1_3} with \eqref{est_k_5} with the previous estimates of this section, we obtain, under our hypotheses on the equation coefficients and the Levi conditions (LC) the energy estimate
\[
\begin{split}
\frac{dE(t)}{dt}&\prec \biggl(E(t)+\Vert g_0\Vert_{H^2}^2+\Vert g_1\Vert_{H^1}^2+\Vert g_2\Vert_{L^2}^2+\int_{0}^t\Vert f(s)\Vert_{L^2}^2\, ds+\int_0^t E(s)\, ds+\Vert f(t)\Vert_{H^1}^2\biggr)\\
&\prec \biggl(E(t)+\int_0^t E(s)\, ds+\Vert f(t)\Vert_{H^1}^2+\Vert f\Vert_{L^\infty\times L^2}^2+\Vert g_0\Vert _{H^2}^2+\Vert g_1\Vert_{H^1}^2+\Vert g_2\Vert_{L^2}^2 \biggr).
\end{split}
\]
This leads to, for $i=5,10$ to 
\[
\begin{split}
\Vert V_{i}\Vert ^2_{L^2}\le E(t)&\prec E(0)+\int_0^t \Vert f(s)\Vert_{H^1}^2\, ds+ \Vert f\Vert_{L^\infty\times L^2}^2+\Vert g_0\Vert _{H^2}^2+\Vert g_1\Vert_{H^1}^2+\Vert g_2\Vert_{L^2}^2\\
&\prec \int_0^t \Vert f(s)\Vert_{H^1}^2\, ds+ \Vert f\Vert_{L^\infty\times L^2}^2+\Vert g_0\Vert^2_{H^3}+\Vert g_1\Vert^2_{H^2}+\Vert g_2\Vert^2_{H^1},
\end{split}
\]
where we used the fact that
\[
E(0)\prec \Vert g_0\Vert _{H^3}^2+\Vert g_1\Vert_{H^2}^2+\Vert g_2\Vert_{H^1}^2.
\]
Hence,
\[
\Vert u(t)\Vert_{H^1}^2\le c\biggl( \int_0^t \Vert f(s)\Vert_{H^1}^2\, ds+ \Vert f\Vert_{L^\infty\times L^2}^2+\Vert g_0\Vert^2_{H^3}+\Vert g_1\Vert^2_{H^2}+\Vert g_2\Vert^2_{H^1}\biggr).
\]
Note that the same kind of inequality holds also for $n>2$. We will simply have to deal with bigger size matrices and estimate carefully the entries $V_k$, with $k=2n+1, 2(2n+1), \cdots, n(2n+1)$.
We have therefore proven the following statement which is an extension of Theorem \ref{theo_3_case_1}(ii) to any space dimension.
\begin{theorem}
\label{theo_n_H_1}
Let us consider the Cauchy problem
\[
\begin{split}
\partial_t^3u-\sum_{i=1}^n a_i(x)\partial_t\partial_{x_i}^2u+\sum_{i=1}^n b_i(x)\partial_{x_i}^2 u+\sum_{i=1}^n  b_{2,i}(x)\partial_t\partial_{x_i}u+b_{3,n}(x)\partial_t^2u&=f(t,x),\\
u(0,x)&=g_0(x),\\
\partial_tu(0,x)&=g_1(x),\\
\partial_t^2u(0,x)&=g_2(x),
\end{split}
\]
where the equation coefficients are real-valued, smooth and with bounded derivatives of any order, $a_i\ge 0$ for $i=1,2$ and $f\in C([0,T],H^1(\R^n))$. Under the Levi conditions (LC),
\[
\begin{split}
b_i&=\lambda a_i,\qquad i=1,\cdots, n,\\
|b_{2,i}|&\prec \sqrt{a_i},\qquad i=1,\cdots, n,\\
|b_{3,n}|&\prec 1,
\end{split}
\]
where $\lambda\in B^\infty(\R^n)$, the Cauchy problem has a unique solution in $C^3([0,T], H^1(\R^n))$ provided that $g_0\in H^3(\R^n)$, $g_1\in H^2(\R^n)$ and $g_2\in H^1(\R^n)$.
\end{theorem}
\begin{remark}
Note that in the proof we have only used the fact that $\lambda$ is bounded with bounded first order derivatives. However, to extend this result to any Sobolev order we will need to assume that $\lambda$ belong to $B^\infty(\R^n)$.
\end{remark}
We now want to prove that the Levi conditions above guarantee well-posedness in every Sobolev space. 
As a first step we start from the system in $V$
\[
\partial_t V=\sum_{i=1}^n \wt{A_i}(x)\partial_{x_i}V+ \wt{B}(x)V+\wt{F}(t,x),
\]
and we derive once more with respect to $x$. We get a system in $W$ where
\[
W=(\partial_{x_1}V,\cdots, \partial_{x_n}V)^T.
\]
$W$ is a $n^2(2n+1)$ column vector with entries $\partial_{x_i}V$, $i=1,\cdots,n$. 
By straitghforward computations we obtain
\[
\partial_{x_j}\partial_t V=\sum_{i=1}^n \wt{A_i}(x)\partial_{x_j}\partial_{x_i}V+\sum_{i=1}^n \partial_{x_j}\wt{A_i}(x)V+\wt{B}(x)\partial_{x_j}V+\partial_{x_j}\wt{B}(x)V+\partial_{x_j}\wt{F}(t,x),
\]
for $j=1,\cdots,n$. In column notation we can therefore write
\[
\partial_t W=\sum_{i=1}^n \wt{\wt{A_i}}(x)\partial_{x_i}W+\wt{\wt{B}}W+\wt{\wt{F}}(t,x),
\]
where, $\wt{\wt{A_i}}$ as a diagonal block structure with repeated block $A_i$,
\[
\wt{\wt{B}}=\left(
	\begin{array}{ccccc}
	\partial_{x_1}{\wt{A_1}}+{\wt{B}} & \partial_{x_1}{\wt{A_2}} & \cdots & \cdots & \partial_{x_1}{\wt{A_n}}\\
	\partial_{x_2}{\wt{A_1}}& \partial_{x_2}{\wt{A_2}}+{\wt{B}} & \cdots & \cdots & \partial_{x_2}{\wt{A_n}}\\
	\vdots & \vdots & \vdots & \vdots & \vdots\\
	\partial_{x_k}{\wt{A_1}}&  \cdots &\partial_{x_k}{\wt{A_k}}+{\wt{B}} & \cdots & \partial_{x_k}{\wt{A_n}}\\
	\vdots & \vdots & \vdots & \vdots & \vdots\\
	\partial_{x_n}{\wt{A_1}}& \partial_{x_n}{\wt{A_2}} & \cdots & \cdots & \partial_{x_n}{\wt{A_n}}+{\wt{B}}\\
	 \end{array}
	\right),
\]
and
\[
{\wt{\wt{F}}}=\nabla_x {\wt{F}}+\left(
	\begin{array}{c}
	(\partial_{x_1}{B})V\\ 
	(\partial_{x_2}{B})V\\ 
	\vdots \\
	(\partial_{x_k}{B})V\\
	\vdots\\
	(\partial_{x_n}{B})V\\ 
	\end{array}
	\right).
\]
This is exactly the same structure we had for the system in $V$ with one more block diagonal step in our argument (denoted with en extra suffix $\wt{\,}$\,) and with $U$ replaced by $V$. We therefore can repeat the same argument employed for the system in $V$ if we make use of the energy $E(t)=(\wt{\wt{Q}}W,W)_{L^2}$. It is also clear that, because of this special structure, our argument can be reiterated as many times we want, obtaining estimates for every Sobolev order $k$. In this specific case, since we have two iterations or in other words we have derived the original system in $U$ twice, we will get estimates in $H^2$ for the solution $u$. More precisely, by arguing as for the system in $V$ we end up with the estimate
\beq
\label{est_u_H_2_n}
\Vert u(t)\Vert_{H^2}^2\le c\biggl( \int_0^t \Vert f(s)\Vert_{H^2}^2\, ds+ \Vert f\Vert_{L^\infty\times H^1}^2+\Vert g_0\Vert^2_{H^4}+\Vert g_1\Vert^2_{H^3}+\Vert g_2\Vert^2_{H^2}\biggr),
\eeq
and, in general, by iteration
\beq
\label{est_u_H_k_n}
\Vert u(t)\Vert_{H^k}^2\le c\biggl( \int_0^t \Vert f(s)\Vert_{H^k}^2\, ds+ \Vert f\Vert_{L^\infty\times H^{k-1}}^2+\Vert g_0\Vert^2_{H^{k+2}}+\Vert g_1\Vert^2_{H^{k+1}}+\Vert g_2\Vert^2_{H^k}\biggr),
\eeq
for all $t\in[0,T]$ and $x\in\R^n$.
\begin{theorem}
\label{theo_n_H_k}
Let us consider the Cauchy problem
\[
\begin{split}
\partial_t^3u-\sum_{i=1}^n a_i(x)\partial_t\partial_{x_i}^2u+\sum_{i=1}^n b_i(x)\partial_{x_i}^2 u+\sum_{i=1}^n  b_{2,i}(x)\partial_t\partial_{x_i}u+b_{3,n}(x)\partial_t^2u&=f(t,x),\\
u(0,x)&=g_0(x),\\
\partial_tu(0,x)&=g_1(x),\\
\partial_t^2u(0,x)&=g_2(x),
\end{split}
\]
where the equation coefficients are real-valued, smooth and with bounded derivatives of any order, $a_i\ge 0$ for $i=1\cdots, n$ and $f\in C([0,T],H^k(\R^n))$. Under the Levi conditions (LC),
\[
\begin{split}
b_i&=\lambda a_i,\qquad i=1,\cdots, n,\\
|b_{2,i}|&\prec \sqrt{a_i},\qquad i=1,\cdots, n,\\
\end{split}
\]
where $\lambda\in B^\infty(\R^n)$, the Cauchy problem has a unique solution in $C^3([0,T], H^k(\R^n))$ provided that $g_0\in H^{k+2}(\R^n)$, $g_1\in H^{k+1}(\R^n)$ and $g_2\in H^k(\R^n)$.
\end{theorem}
It clearly follows that the Cauchy problem above is $C^\infty$ well-posed.


In this paper we have focused on a specific third order equation to explain better how to overcome the technical difficulties coming with a higher space dimension. We plan to address the general order $m>3$ in a forthcoming paper. We expect our method to still hold independently of the equation order. 

\section{Examples of higher order hyperbolic equations in any space dimension}

We conclude this paper by testing our method on some higher order examples related to the equation studied in the previous section. 

We want to study the well-posedness of the Cauchy problem for equations of the type
\begin{multline*}
\partial_t^4u-\sum_{i=1}^n a_i(x)\partial_t^2\partial_{x_i}^2u+\sum_{i=1}^n b_i(x)\partial_{x_i}^3 u+\sum_{i=1}^n  b_{2,i}(x)\partial_t\partial^2_{x_i}u\\
+\sum_{i=1}^n b_{3,i}(x)\partial_t^2\partial_{x_i}u+b_{4,n}(x)\partial_t^3u=f(t,x),
\end{multline*}
where $a_i\ge 0$ for all $i=1,\cdots,n$. Note that in one space dimension the equation above would be associated to a system of differential equations with Sylvester matrix 
\[
A=\left(\begin{array}{cccc}
0 & 1 & 0 & 0 \\
0& 0 & 1 & 0\\
0 & 0 & 0 & 1\\
0 & 0 & a & 0
 \end{array}
	\right)
\]
and symmetriser
\[
Q=\left(\begin{array}{cccc}
0 & 0 & 0 & 0 \\
0& 2a^2 & 0 & -2a\\
0 & 0 & 2a & 0\\
0 & -2a & 0 & 4
 \end{array}
	\right).
\]
By direct computations,
\[
QA=\left(\begin{array}{cccc}
0 & 0 & 0 & 0 \\
0& 0 & 0 & 0\\
0 & 0 & 0 & 2a\\
0 & 0 & 2a & 0
 \end{array}
	\right).
\]
In higher dimensions, the equation 
\[
\partial_t^4-\sum_{i=1}^n a_i(x)\partial_t^2\partial_{x_i}^2+\sum_{i=1}^n b_i(x)\partial_{x_i}^3 u+\sum_{i=1}^n  b_{2,i}(x)\partial_t\partial^2_{x_i}u+\sum_{i=1}^n b_{3,i}(x)\partial_t^2\partial_{x_i}u+b_{4,n}(x)\partial_t^3u=f(t,x),
\]
is transformed via
\[
U=(\partial^3_{x_1}u,\cdots, \partial^3_{x_n}u, \partial_t\partial^2_{x_1}u,\cdots, \partial_t\partial^2_{x_n}u, \partial_t^2\partial_{x_1}u,\cdots,\partial_t^2\partial_{x_n}u, \partial^3_tu)^T
\]
into the $3n+1\times 3n+1$ system
\[
\partial_t U=\sum_{k=1}^n A_k(x)\partial_{x_k}U+B(x)U+F,
\]
where, the matrices $A_k$ have entries $(a_{k,ij})_{ij}$ as follows:
\[
\begin{split}
a_{k,ij}&=1,\quad\text{for $i=k$ and $j=n+k$},\\
a_{k,ij}&=1,\quad\text{for $i=n+k$ and $j=2n+k$},\\
a_{k,ij}&=1,\quad\text{for $i=2n+k$ and $j=3n+1$},\\
a_{k,ij}&=a_k,\quad\text{for $i=3n+1$ and $j=2n+k$},\\
a_{k,ij}&=0,\quad\text{otherwise},
\end{split}
\]
$B$ has all the rows vanishing a part from the last given by 
\[
(-b_1,\cdots,-b_n, -b_{2,1},\cdots, -b_{2,n}, -b_{3,1},\cdots,-b_{3,n}, -b_{4,n})
\]
and $F$ is the column vector with entries $f_{i,1}=0$ if $i\neq 3n=1$ and $f_{3n+1,1}=f$.

Without loss of generality, we assume that $n=2$. We therefore, deal with matrices ${A}_i$ in Sylvester form where the part in bold coincides with the Sylvester matrices encountered in the previous section, i.e., 
\[
{A}_1=\left(\begin{array}{ccccccc}
0 & 0 & 1 & 0 & 0 & 0 & 0\\
0 & 0 & 0 & 0 & 0 & 0 & 0\\
0 & 0 & {\bf 0} & {\bf 0} & {\bf 1} & {\bf 0} & {\bf 0}\\
0 & 0 & {\bf 0} & {\bf 0} & {\bf 0}& {\bf 0}& {\bf 0}\\ 
0 & 0 & {\bf 0} & {\bf 0} & {\bf 0}& {\bf 0} & {\bf 1}\\
0 & 0 & {\bf 0} & {\bf 0} & {\bf 0}& {\bf 0} & {\bf 0}\\
0 & 0 & {\bf 0} & {\bf 0}& {\bf a_1} & {\bf 0} & {\bf 0}\\
	\end{array}
	\right), 
{A}_2=\left(\begin{array}{ccccccc}
0 & 0 & 0 & 0 & 0 & 0 & 0\\
0 & 0 & 0 & 1 & 0 & 0 & 0\\
0 & 0 & {\bf 0} & {\bf 0} & {\bf 0} & {\bf 0} & {\bf 0}\\
0 & 0 & {\bf 0} & {\bf 0} & {\bf 0}& {\bf 1}& {\bf 0}\\ 
0 & 0 & {\bf 0} & {\bf 0} & {\bf 0}& {\bf 0} & {\bf 0}\\
0 & 0 & {\bf 0} & {\bf 0} & {\bf 0}& {\bf 0} & {\bf 1}\\
0 & 0 & {\bf 0} & {\bf 0}& {\bf 0} & {\bf a_2} & {\bf 0}\\
	\end{array}
	\right).
\]
Inspired by the method applied to third order equations, we define 
\[
Q=\frac{1}{3}\left(\begin{array}{ccccccc}
0 & 0 & 0 & 0 & 0 & 0 & 0\\
0 & 0 & 0 & 0 & 0 & 0 & 0\\
0 & 0 & a_1^2& a_1a_2 & 0 &  0& -a_1\\
0 & 0 & a_1a_2& a_2^2 & 0& 0& -a_2\\ 
0 & 0 & 0 & 0 & 2a_1& 0 & 0\\
0 & 0 & 0 & 0 & 0 & 2a_2& 0\\
0 & 0 & -a_1 & -a_2& 0 & 0 & 3\\
	\end{array}
	\right).
\]
We immediately have
\[
QA_1=\left(\begin{array}{ccccccc}
0 & 0 & 0 & 0 & 0 & 0 & 0\\
0 & 0 & 0 & 0 & 0 & 0 & 0\\
0 & 0 & 0 & 0 & 0 & 0 & 0\\
0 & 0 & 0 & 0 & 0 & 0 & 0\\ 
0 & 0 & 0 & 0 &  0& 0 & 2a_1\\
0 & 0 & 0 & 0 & 0 & 0 & 0\\
0 & 0 & 0 & 0& 2a_1 & 0 & 0\\
	\end{array}
	\right), QA_2=\left(\begin{array}{ccccccc}
0 & 0 & 0 & 0 & 0 & 0 & 0\\
0 & 0 & 0 & 0 & 0 & 0 & 0\\
0 & 0 & 0 & 0 & 0 & 0 & 0\\
0 & 0 & 0 & 0 & 0 & 0 & 0\\ 
0 & 0 & 0 & 0 &  0& 0 & 0\\
0 & 0 & 0 & 0 & 0 & 0 & 2a_2\\
0 & 0 & 0 & 0& 0 & 2a_2 & 0\\
	\end{array}
	\right).
\]
In general, we have the following proposition which follows immediately from Proposition \ref{prop_symm_3_n} with a $n$-shift in the indexes. 
\begin{proposition}
\label{prop_symm_4_n}
Let  $Q$ the block diagonal $3n+1\times 3n+1$ matrix defined by a $n\times n$ zero block and the $2n+1\times 2n+1$ block
\[
Q_{2n+1}=\frac{1}{3}\left(
	\begin{array}{cccccccc}
	a_1^2 & a_1a_2& \cdots &a_1a_{n} & 0& \cdots & \cdots & -a_1\\
	a_1a_2 & a_2^2 & a_2a_3 & \cdots & \cdots &  0 & \cdots & -a_2\\
	\vdots & \vdots & \vdots & \vdots & \vdots & \vdots & \vdots & \vdots\\
	a_1a_k & \vdots & a_k^2  & \vdots& 0 & \vdots  &\vdots& -a_k\\
	\vdots & \vdots & \vdots & \vdots & \vdots & \vdots & \vdots & \vdots\\
	a_1a_n & \vdots & \vdots &a_n^2 &  \vdots &\vdots & \vdots & -a_n\\
	\vdots & \vdots & \vdots & \vdots & 2a_1 & \vdots & \vdots & 0\\
	\vdots & \vdots & \vdots & \vdots & \vdots & \vdots & \vdots & \vdots\\
	\vdots & \vdots & \vdots & \vdots & \vdots & \vdots & 2a_n &  0\\
	-a_1 & -a_2 & \cdots & -a_n & 0 & \vdots & 0 & 3
	 \end{array}
	\right),
\]
i.e., 
\[
Q=\left(
	\begin{array}{cc}
	0 & 0\\
	0 & Q_{2n+1}
	\end{array}
	\right).
\]
\begin{itemize}
\item[(i)] $Q$ is a symmetriser of $A_k$ for every $k=1,\dots,n$, i.e, $QA_k=A_k^\ast Q$ is the symmetric matrix with entries $2n+k, 3n+1$ and $3n+1, 2n+k$ equals to $\frac{2}{3}a_k$ and $0$ otherwise.
\item[(ii)] For every $v\in \R^{3n+1}$
\[
\begin{split}
3\lara{Qv,v}&=\sum_{k=1}^n\lara{a_k^2v_{n+k},v_{n+k}}+2\sum_{k=1}^n \lara{a_k v_{2n+k},v_{2n+k}}+3\Vert v_{3n+1}\Vert^2-2\sum_{k=1}^n\lara{a_kv_{n+k}, v_{3n+1}}\\
&+2\sum_{1\le i<j\le n}\lara{a_iv_{n+i},a_{j}v_{n+j}}\\
&=\Vert \sum_{k=1}^n a_{n+k}v_k-v_{3n+1}\Vert^2+2\sum_{k=1}^n \lara{a_k v_{2n+k},v_{2n+k}}+2\Vert v_{3n+1}\Vert^2.
\end{split}
\]
\end{itemize}
\end{proposition}
 The block-diagonal structure of the symmetriser reduces the analysis of the Cauchy problem 
\beq
\label{4_n_eq}
\begin{split}
\partial_t^4u-\sum_{i=1}^n a_i(x)\partial_t^2\partial_{x_i}^2u&+\sum_{i=1}^n b_i(x)\partial_{x_i}^3 u+\sum_{i=1}^n  b_{2,i}(x)\partial_t\partial^2_{x_i}u+\sum_{i=1}^n b_{3,i}(x)\partial_t^2\partial_{x_i}u+b_{4,n}(x)\partial_t^3u=f(t,x),\\
u(0,x)&=g_0(x),\\
\partial_tu(0,x)&=g_1(x),\\
\partial_t^2u(0,x)&=g_2(x),\\
\partial_t^3u(0,x)&=g_3(x),
\end{split}
\eeq
or equivalently of the Cauchy problem
\[
\begin{split}
\partial_t U&=\sum_{k=1}^n A_k(x)\partial_{x_k}U+B(x)U+F,\\
U(0)&=(\partial^3_{x_1}g_0,\cdots, \partial^3_{x_n}g_0, \partial_{x_1}^2g_1,\cdots, \partial_{x_n}^2g_1, \partial_{x_1}g_2,\cdots,\partial_{x_n}g_2, g_3)^T\\
\end{split}
\]
to the analysis of the Cauchy problem \eqref{CP_3_n} in the previous section.  The Levi conditions on the lower order terms are obtained by setting 
\[
3((QB+B^\ast Q)U,U)_{L^2}\prec E(t)=(QU, U)_{L^2}.
\]
Without loss of generality, we can assume that $n=2$. We have 
\[
3QB=\left(
\begin{array}{cccccccc}
0 & 0 & \cdots & \cdots & \cdots & 0 & 0\\
0 & 0 & \cdots & \cdots & \cdots & 0 & 0\\
a_1b_1 & a_1b_2 & a_1 b_{2,1} & a_1b_{2,2} & a_1 b_{3,1} & a_1 b_{3,2} & a_1b_{4,2}\\
a_2b_1 & a_2b_2 & a_2 b_{2,1} & a_2b_{2,2} & a_2 b_{3,1} & a_2 b_{3,2} & a_2b_{4,2}\\
0 & 0 & \cdots & \cdots & \cdots & 0 & 0\\
0 & 0 & \cdots & \cdots & \cdots & 0 & 0\\
-3b_1 & -3b_2 & -3b_{2,1} & -3b_{2,2} & -3b_{3,1} & -3b_{3,2} & -3b_{4,2}\
 \end{array}
	\right)
\]
and $3(QB+B^\ast Q)$ is the matrix
\[
\left(
\begin{array}{cccccccc}
0 & 0 & a_1b_1 & a_2b_1 & 0 & 0 & -3b_1\\
0 & 0 & a_1b_2 & a_2b_2 & 0 & 0 & -3b_2\\
a_1b_1 & a_1b_2 & 2a_1 b_{2,1} & a_1b_{2,2}+a_2b_{2,1} & a_1 b_{3,1} & a_1 b_{3,2} & -3b_{2,1}+a_1b_{4,2}\\
a_2b_1 & a_2b_2 & a_1b_{2,2}+a_2 b_{2,1} & 2a_2b_{2,2} & a_2 b_{3,1} & a_2 b_{3,2} & -3b_{2,2}+a_2b_{4,2}\\
0 & 0 & a_1b_{3,1} & a_2b_{3,1} & 0 & 0 & -3b_{3,1}\\
0 & 0 & a_1b_{3,2} & a_2b_{3,2} & 0 & 0 & -3b_{3,2}\\
-3b_1 & -3b_2 & -3b_{2,1}+a_1b_{4,2} & -3b_{2,2}+a_2b_{4,2} & -3b_{3,1} & -3b_{3,2} & -6b_{4,2}\
 \end{array}
	\right).
\]
It follows that 
\[
\begin{split}
3((QB+B^\ast Q)U,U)_{L^2}
&= 2(a_1b_1U_3,U_1)_{L^2}+2(a_2b_1U_4,U_1)_{L^2}+2(-3b_1U_7,U_1)_{L^2}\\
&+2(a_1b_2U_3,U_2)_{L^2}+2(a_2b_2U_4,U_2)_{L^2}+2(-3b_2U_7,U_2)_{L^2}\\
&+(2a_1b_{2,1}U_3,U_3)_{L^2}+2((a_1b_{2,2}+a_2 b_{2,1})U_4,U_3)_{L^2}+2(a_1b_{3,1}U_5,U_3)_{L^2}\\
&+2(a_1b_{3,2}U_6,U_3)_{L^2}+2((-3b_{2,1}+a_1b_{4,2})U_7,U_3)_{L^2}\\
&+(2a_2b_{2,2}U_4,U_4)_{L^2}+2(a_2b_{3,1}U_5,U_4)_{L^2}\\
&+2(a_2b_{3,2}U_6,U_4)_{L^2}+2((-3b_{2,2}+a_2b_{4,2})U_7,U_4)_{L^2}\\
&+2(-3b_{3,1}U_7,U_5)_{L^2}+2(-3b_{3,2}U_7,U_6)_{L^2}+(-6b_{4,2}U_7,U_7)_{L^2}.
\end{split}
\]
Note that 
\[
\begin{split}
3E(t)&=(a_1^2U_3,U_3)_{L^2}+(a_1a_2U_4,U_3)_{L^2}-2(a_1U_7,U_3)_{L^2}+(a_2^2U_4,U_4)_{L^2}-2(a_2U_7,U_4)_{L^2}\\
&+2(a_1U_5,U_5)_{L^2}+2(a_2U_6,U_6)_{L^2}+3\Vert U_7\Vert_{L^2}^2\\
&= \Vert a_1U_3+ a_2U_4-U_7\Vert_{L^2}^2+2(a_1U_5,U_5)_{L^2}+2(a_2U_6,U_6)_{L^2}+2\Vert U_7\Vert_{L^2}^2\\
&\ge \frac{1}{2}\Vert a_1U_3+ a_2U_4\Vert_{L^2}^2+2(a_1U_5,U_5)_{L^2}+2(a_2U_6,U_6)_{L^2}+\Vert U_7\Vert_{L^2}^2.
\end{split}
\]
By setting 
\[
((QB+B^\ast Q)U,U)_{L^2}\prec E(t)
\]
we obtain that since $U_1$ and $U_2$ do not appear in $E(t)$ necessarily $b_1=b_2=0$. We now want to find suitable Levi conditions which allow to estimate $3((QB+B^\ast Q)U,U)_{L^2}$ with 
\[
\frac{1}{2}\Vert a_1U_3+ a_2U_4\Vert_{L^2}^2+2(a_1U_5,U_5)_{L^2}+2(a_2U_6,U_6)_{L^2}+\Vert U_7\Vert_{L^2}^2.
\]
We rewrite $3((QB+B^\ast Q)U,U)_{L^2}$ as the sum of 
\[
\begin{split}
&I_1=(2a_1b_{2,1}U_3,U_3)_{L^2}+2((a_1b_{2,2}+a_2 b_{2,1})U_4,U_3)_{L^2}+(2a_2b_{2,2}U_4,U_4)_{L^2}\\
&I_2=2(a_1b_{3,1}U_5,U_3)_{L^2}+2(a_2b_{3,1}U_5,U_4)_{L^2}\\\
&I_3=2(a_1b_{3,2}U_6,U_3)_{L^2}+2(a_2b_{3,2}U_6,U_4)_{L^2}\\
&I_4=2((-3b_{2,1}+a_1b_{4,2})U_7,U_3)_{L^2}+2((-3b_{2,2}+a_2b_{4,2})U_7,U_4)_{L^2}\\
&I_5=2(-3b_{3,1}U_7,U_5)_{L^2}+2(-3b_{3,2}U_7,U_6)_{L^2}+(-6b_{4,2}U_7,U_7)_{L^2}.
\end{split}
\]
By setting 
\[
\begin{split}
b_{2,2}&=\lambda a_2,\\
b_{2,1}&=\lambda a_1,
\end{split}
\]
we easily see that 
\[
I_1\prec \Vert |\lambda|a_1U_3+|\lambda|a_2U_4\Vert_{L^2}^2\prec E(t).
\]
Since 
\[
I_2=2(b_{3,1}U_5,a_1U_3+a_2U_4)_{L^2}\prec \Vert a_1U_3+a_2U_4\Vert_{L^2}^2+(b_{3,1}^2U_5,U_5)_{L^2}
\]
by imposing 
\[
|b_{3,1}|\prec \sqrt{a_1}
\]
we have that 
\[
I_2\prec E(t).
\]
Analogously,
\[
I_3=2(b_{3,2}U_6,a_1U_3+a_2U_4)_{L^2} 
\]
leads to the Levi condition
\[
|b_{3,2}|\prec \sqrt{a_2}.
\]
We now write $I_4$ as 
\[
\begin{split}
&2((-3\lambda a_1+a_1b_{4,2})U_7,U_3)_{L^2}+2((-3\lambda a_2+a_2b_{4,2})U_7,U_4)_{L^2}\\
&2(-3\lambda U_7, a_1U_3+a_2U_4)_{L^2}+2(2b_{4,2}U_7, a_1U_3+a_2U_4)_{L^2}.
\end{split}
\]
Hence, if $|b_{4,2}|\prec 1$ we obtain that 
\[
I_4\prec \Vert U_7\Vert_{L^2}^2+\Vert a_1U_3+a_2U_4\Vert_{L^2}^2\prec E(t).
\]
Finally, analysing $I_5$ we easily see that the Levi conditions $b_{3,1}\prec \sqrt{a_1}$, $b_{3,2}\prec \sqrt{a_2}$ and $|b_{4,2}|\prec 1$ imply $I_5\prec E(t)$. Concluding, our method allow us to identify the following Levi conditions for the equation \eqref{4_n_eq}:
\[
\begin{split}
b_i&=0,\\
b_{2,i}&=\lambda a_i,\\
|b_{3,i}|&\prec \sqrt{a_i},\\
|b_{4,n}|&\prec 1
\end{split}
\]
for all $i=1,\cdots, n$, with $\lambda\in B^\infty(\R^n)$.
 
Note that this is an extension to higher space dimension of the Levi conditions formulated in \cite{ST21} in dimension 1. Indeed, the equation 
\[
\partial_t^4u-a(x)\partial_t^2\partial_{x}^2u+b_1(x)\partial_{x}^3 u+b_{2}(x)\partial_t\partial^2_{x}u+b_{3}(x)\partial_t^2\partial_{x}u+b_{4}(x)\partial_t^3u=f(t,x)
\]
has roots 
\[
-\sqrt{a},\, 0,\, 0\, \sqrt{a}
\]
and lower order terms associated to the polynomial
\[
R(\tau,x)=b_1(x)+b_2(x)\tau+b_3(x)\tau^2+b_4(x)\tau^3.
\]
The Levi conditions formulated in \cite{ST21} in order to get $C^\infty$ well-posedness requires that 
\[
\begin{split}
&b_1(x)+b_2(x)\tau+b_3(x)\tau^2+b_4(x)\tau^3\\
&=l_1(x)\tau^2(\tau-\sqrt{a(x)})+(l_2(x)+l_3(x))\tau(\tau^2-a(x))+l_4(x)\tau^2(\tau+\sqrt{a(x)}),
\end{split}
\]
where the functions $l_i$ are bounded. Hence
\[
\begin{split}
&b_1(x)+b_2(x)\tau+b_3(x)\tau^2+b_4(x)\tau^3\\
&=-(l_2+l_3)(x)a(x)\tau+(-l_1+l_4)(x)\sqrt{a(x)}\tau^2+(l_1+l_2+l_3+l_4)(x)\tau^3.
\end{split}
\]
It follows that 
\[
\begin{split}
b_1=&0,\\
b_2(x)=& -(l_2+l_3)(x)a(x),\\
b_3(x)=&(-l_1+l_4)(x)\sqrt{a(x)},\\
|b_4|&\prec 1.
\end{split}
\]
These are a special case of our Levi conditions formulated in $\R^n$.

\subsection{Conclusion}
The Levi conditions deduced in this paper force the equation in \eqref{4_n_eq} to be written in the simpler form 
\[
\partial_t^4u-\sum_{i=1}^n a_i(x)\partial_t^2\partial_{x_i}^2u+\sum_{i=1}^n  b_{2,i}(x)\partial_t\partial^2_{x_i}u+\sum_{i=1}^n b_{3,i}(x)\partial_t^2\partial_{x_i}u+b_{4,n}(x)\partial_t^3u=f(t,x).
\]
The study of the corresponding Cauchy problem can be therefore reduced to the third order model investigated in the previous section. We immediately obtain the following well-posedness result.
\begin{theorem}
\label{theo_4_final}
Let 
\[
\begin{split}
\partial_t^4u-\sum_{i=1}^n a_i(x)\partial_t^2\partial_{x_i}^2u&+\sum_{i=1}^n b_i(x)\partial_{x_i}^3 u+\sum_{i=1}^n  b_{2,i}(x)\partial_t\partial^2_{x_i}u+\sum_{i=1}^n b_{3,i}(x)\partial_t^2\partial_{x_i}u+b_{4,n}(x)\partial_t^3u=f(t,x),\\
u(0,x)&=g_0(x),\\
\partial_tu(0,x)&=g_1(x),\\
\partial_t^2u(0,x)&=g_2(x),\\
\partial_t^3u(0,x)&=g_3(x),
\end{split}
\]
where the equation coefficients are real-valued, smooth and with bounded derivatives of any order, $a_i\ge 0$ for $i=1\cdots, n$ and $f\in C([0,T],H^k(\R^n))$. Under the Levi conditions (LC),
\[
\begin{split}
b_i&=0,\\
b_{2,i}&=\lambda a_i,\quad i=1,\cdots,n\\
|b_{3,i}|&\prec \sqrt{a_i},\quad i=1,\cdots,n\\
\end{split}
\]
where $\lambda\in B^\infty(\R^n)$, the Cauchy problem has a unique solution in $C^3([0,T], H^k(\R^n))$ provided that $g_0\in H^{k+3}(\R^n)$, $g_1\in H^{k+2}(\R^n)$, $g_2\in H^{k+1}(\R^n)$ and $g_3\in H^{k}(\R^n)$.
\end{theorem}
\begin{remark}
Note that the same kind of Levi conditions appear also in higher order equations with coefficients in $B^\infty(\R^n)$ which can be easily reduced to the third order model. For instance, for $m\in\N$, let us consider the equation
\[
\partial^{2m}_t u-\sum_{i=1}^n a_i(x)\partial^{2m-2}_t\partial_{x_i}^2u=f(t,x),
\]
where $a_i\ge 0$ for all $i=1,\cdots,n$. If we now add lower order terms of the type 
\[
\sum_{i=1}^n b_{i}(x)\partial_t^{2m-3}\partial_{x_i}^2u+\sum_{i=1}^n b_{2,i}(x)\partial_t^{2m-2}\partial_{x_i}u+b_{3,n}(x)\partial_t^{2m-1}u
\]
then the corresponding Cauchy problem is $C^\infty$ well-posed provided that 
\[
\begin{split}
b_i&=\lambda a_i,\qquad i=1,\cdots, n,\\
|b_{2,i}|&\prec \sqrt{a_i},\qquad i=1,\cdots, n,\\
\end{split}
\]
where $\lambda\in B^\infty(\R^n)$.
\end{remark}

\section{Appendix: symmetriser of a matrix in Sylvester form}
In this sequel we collect some important properties about standard symmetrisers. These properties are mainly proven in \cite{JT} so we do not give detailed proofs. However, a sketch of a proof is provided for those statements which are of fundamental importance for the paper. Note that for the physical meaning of our problem we assume that all the functions and matrices are  {\bf real-valued}. We start with the following general algebraic results valid for positive definite and semidefinite matrices. Let $\lara{\cdot,\cdot}$ be the Euclidean scalar product in $\R^n$. We recall that an $m\times m$ symmetric matrix $M$ is positive definite if $\lara{Mv,v}>0$ for all $v\in\R^m$ and positive semidefinite if $\lara{Mv,v}\ge 0$ for all $v\in\R^m$. Positive definite and semidefinite matrices are characterised by their eigenvalues, in the sense that $A$ is positive (semi)definite if and only if all its eigenvalues are (non-negative) positive. A matrix $A$ is in Sylvester form if it is of the type
\[
A=\left(
	\begin{array}{cccc}
	0& 1 & \cdots & 0\\
	0& 0 & 1 & \cdots\\
	\vdots & \vdots & \vdots & \vdots\\
        	a_{m1} & a_{m2} & \cdots & a_{mm}
	\end{array}
	\right).
\]
\begin{proposition}
\label{prop_exist_sym}
\leavevmode
\begin{itemize}
\item[(i)] Every Sylvester matrix $A$ with real eigenvalues $\lambda_1, \lambda_2,\dots,\lambda_m$ admits a symmetriser, i.e., a symmetric matrix $Q$ such that $AQ=QA^\ast$. The entries of $Q$ are polynomials in the eigenvalues $\lambda_i$, $i=1,\dots,m$ and can also be written as polynomials in $a_{m1}, a_{m2},\dots, a_{mm}$.
\item[(ii)] Let $\sigma_{0,k}=1$ and 
\[
\sigma_{h,k}(\lambda):=(-1)^h\sum_{\substack{1\le j_1<\cdots<j_h\le m\\ j_i\neq k}}\lambda_{j_1}\cdots\lambda_{j_h},
\]
with $1\le h,k\le m$. Hence, the symmsetriser $Q(\lambda)$ in $(ii)$ has entries
\[
q_{ij}=m^{-1}\sum_{1\le k\le m}\sigma_{m-i,k}(\lambda)\sigma_{m-j,k}(\lambda),\qquad 1\le i,j\le m.
\]
\item[(iii)] If
\[
\Delta:=\prod _{1\le i,j\le m}(\lambda_i-\lambda_j)^2>0
\]
then $Q$ is positive definite. If $\Delta=0$ then $Q$ is positive semi-definite. Moreover, $\det Q=m^{-m}\Delta$.
\item[(iv)] Let
\[
\psi_k(\lambda)=\sum_{1\le j_1<j_2<\cdots<j_k\le m}\lambda_{j_1}^2\lambda_{j_2}^2\cdots \lambda_{j_k}^2,
\]
for $k=1,\dots,m$. Let set $\psi_0=1$ and let $\Psi$ be the $m\times m$ diagonal matrix with entries $\psi_{ii}=\psi_{m-i}$, $i=1,\dots, m$. Hence, there exists a constant $C_m>0$ depending only on the matrix size $m$ such that
\[
\lara{Qv,v}\le C_m\lara{\Psi v,v},
\]
for all $v\in\R^m$.
\end{itemize}
\end{proposition}
\begin{proof}
For assertions (i)-(iii) we refer the reader to \cite{JT} were the symmetriser and its properties are presented and discussed in details. We  observe that by definition of the entries of $Q$ and $\Psi$ there exists a constant $c_m>0$ (depending only on the matrix size $m$) such that 
\[
|\sigma_{m-i,k}(\lambda)|\le c_m\sqrt {\psi_{m-i}(\lambda)},
\]
for all $1\le i,k\le m$. From the second assertion of this proposition it follows that 
\[
|q_{ij}|\le m^{-1}\sum_{1\le k\le m}|\sigma_{m-i,k}(\lambda)||\sigma_{m-j,k}(\lambda)|\le c_m^2\sqrt {\psi_{m-i}(\lambda)}\sqrt {\psi_{m-j}(\lambda)}, \quad 1\le i,j\le m.
\]
Hence,
\[
\begin{split}
\lara{Qv,v}=\sum_{i=1}^m\sum_{j=1}^m q_{ij}v_jv_i&\le \sum_{i=1}^m\sum_{j=1}^m |q_{ij}||v_j||v_i|\\
&\le c_m^2 \sum_{i=1}^m\sum_{j=1}^m \sqrt {\psi_{m-i}(\lambda)}\sqrt {\psi_{m-j}(\lambda)}|v_j||v_i|\\
&= c_m^2\sum_{i=1}^m \psi_{m-i}(\lambda)v_i^2+c_m^2 \sum_{i=1}^m\sum_{j=1, j\neq i}^m \sqrt {\psi_{m-i}(\lambda)}\sqrt {\psi_{m-j}(\lambda)}|v_j||v_i|\\
&\le 2c_m^2\sum_{i=1}^m \psi_{m-i}(\lambda)v_i^2=2c_m^2\lara{\Psi v,v},
\end{split}
\]
for all $v\in\R^m$. This proves (iv) with $C_m=2c_m^2$.
\end{proof}

We now work under the assumption that the eigenvalues $\lambda_i$, $i=1,\dots,m$, fulfil the following property introduced by Kinoshita and Spagnolo in \cite{KS}: there exists $M>0$ such that 
\beq
\label{KS_cond}
\lambda_i^2+\lambda_j^2\le M(\lambda_i-\lambda_j)^2,
\eeq
for all $1\le i<j\le m$. As proven in \cite{JT} this is equivalent to each of the following properties:
\beq
\label{KS1}
\exists M_1>0: \quad \prod_{1\le i<j\le m}(\lambda_i^2+\lambda_j^2)\le M_1\Delta(\lambda),
\eeq
\beq
\label{KS2}
\exists M_2>0: \quad \psi_1(\lambda)\psi_2(\lambda)\cdots \psi_{m-1}(\lambda)\le M_2\Delta(\lambda),
\eeq
for all $\lambda=(\lambda_1,\dots,\lambda_m)\in\R^m$.
 We now recall an important properties of positive definite matrices which will be applied to $Q$ and $\Psi$.
\begin{proposition}
\label{prop_sym_bound}
Let $A_1$ and $A_2$ be two positive definite, symmetric $m\times m$ matrices of order m such
that, for some $c_1, c_2>0$,
\[
\lara{A_1v,v}\le c_1\lara{A_2v,v}\quad \text{for all $v\in\R^m$}
\]
and
\[
{\rm det}A_1\ge c_2{\rm det}A_2.
\]
Then, 
\[
\lara{A_1v,v}\ge c_1^{1-m}c_2\lara{A_2v,v},
\]
for all $v\in \R^m$. 
\end{proposition}
Let us now set $A_1=Q$ and $A_2=\Psi$. Then,
\begin{itemize}
\item there exists a constant $C_1=C_1(m)>0$ such that $\lara{Qv,v}\le C_1\lara{\Psi v,v}$;
\item under the hypothesis \eqref{KS_cond} there exists a constant $C_2>0$ such that 
\[
\det Q =m^{-m}\Delta\ge m^{-m}M_2(\psi_1(\lambda)\psi_2(\lambda)\cdots \psi_{m-1})= C_2\det \Psi;
\]
\item $Q$ and $\Psi$ are positive definite when $\Delta(\lambda)\neq 0$.
\end{itemize}
It follows that we can apply Proposition \ref{prop_sym_bound} directly to $Q=Q(\lambda)$ and $\Psi=\Psi(\lambda)$ when the condition \eqref{KS_cond} holds and $\Delta(\lambda)\neq 0$. Thus,
\[
\lara{Q(\lambda)v,v}\ge C_1^{1-m}C_2\lara{\Psi(\lambda)v,v},
\]
when $\Delta(\lambda)\neq 0$. When $\Delta(\lambda)=0$ we proceed with an approximation argument used already in \cite{KS}. 
Indeed if $\Delta(\lambda)=0$ there exists $k=1,\dots ,m$ such that (up to rearrangement) $\lambda_1=\lambda_2=\dots=\lambda_{k-1}=0$ and $0<|\lambda_k|<|\lambda_{k+1}|<\dots<|\lambda_m|$. We now approximate $\lambda_i$ with $\lambda_{i,\eps}=\eps \frac{\lambda_k}{k}i$ for $i=1,\dots,k-1$ and $\lambda_{i,\eps}=\lambda_i$ for $i=k,\dots,m$. By direct computations we see that $\Delta(\lambda_\eps)\neq 0$ and that \eqref{KS_cond} holds with some constant $M'>0$ independent of $\eps$. Hence, by Proposition \ref{prop_sym_bound}  there exist $C_1, C_2'>0$ such that 
\[
\lara{Q(\lambda_\eps)v,v}\ge C_1^{1-m}C'_2\lara{\Psi(\lambda_\eps)v,v}
\]
for all $\eps\in(0,1]$. By a continuity argument we can therefore conclude that 
\[
\lara{Q(\lambda)v,v}\ge C_1^{1-m}C'_2\lara{\Psi(\lambda)v,v},
\]
for all $v\in\R^m$. We summarise the result proven above in the following proposition.
\begin{proposition}
\label{prop_bound_qpsi}
Under condition \eqref{KS_cond} there exist constants $\gamma_1,\gamma_2>0$ such that
\[
\gamma_1\lara{\Psi(\lambda) v,v}\le \lara{Q(\lambda)v,v}\le \gamma_2\lara{\Psi(\lambda)v,v},
\]
for all $v\in\R^m$.
\end{proposition}
\begin{remark}
\label{rem_imp}
A careful inspection of the argument given above shows that the constants $\gamma_1$ and $\gamma_2$ depend only on the size of the matrices involved and on the constant appearing in \eqref{KS_cond} and its equivalent forms. So, if $\lambda$ is a function of space and time, i.e., $\lambda=\lambda(t,x)$ (or, analogously, if it depends on some parameter), the results above will still hold as long as \eqref{KS_cond} does uniformly in all the variables. More precisely, if there exists $M>0$ such that 
\beq
\label{KS_cond_tx}
\lambda_i^2(t,x)+\lambda_j^2(t,x)\le M(\lambda_i(t,x)-\lambda_j(t,x))^2,
\eeq
\end{remark} 
holds for all $t\in[0,T]$ and $x\in\R^n$ then there exist $\gamma_1,\gamma_2>0$ such that
\beq
\label{prop_est_tx}
\gamma_1\lara{\Psi(\lambda(t,x)) v,v}\le \lara{Q(\lambda(t,x))v,v}\le \gamma_2\lara{\Psi(\lambda(t,x))v,v},
\eeq
for all $v\in\R^m$, uniformly in $t\in[0,T]$ and $x\in\R^n$.

\end{document}